\documentclass[11pt]{article}
\usepackage{latexsym}
\usepackage{graphicx,amsmath,amssymb}
\usepackage{latexsym,array}
\usepackage{multirow}
\usepackage[hang,small,bf]{caption}
\usepackage{mathabx}
\usepackage{float}
\restylefloat{table}

\newtheorem{thm}{Theorem}
\newtheorem{rmk}[thm]{Remark}

\newcommand{\jmod}[2]{\mbox{ $\equiv #1$}\mbox{ \rm (mod $#2$)}}

\date{}

\title{Self-orthogonal codes constructed from weakly self-orthogonal designs invariant under an action of $M_{11}$}

\author{ Vedrana Mikuli{\'c} Crnkovi{\'c} {and} Ivona Traunkar}  

\begin{document}
\maketitle

\begin{abstract}
In this paper we generalize the construction of binary self-orthogonal codes obtained from weakly self-orthogonal designs described by Tonchev in \cite{0} in order to obtain self-orthogonal codes over an arbitrary field. We extend construction self-orthogonal codes from orbit matrices of self-orthogonal designs and weakly self-orthogonal 1-designs such that block size is odd and block intersection numbers are even described in \cite{2}. Also, we generalize mentioned construction in order to obtain self-orthogonal codes over an arbitrary field. We construct weakly self-orthogonal designs invariant under an action of Mathieu group $M_{11}$ and, from them, binary self-orthogonal codes.
\end{abstract}

\section{Introduction}
Codes constructed from block designs (\cite{AKey}, \cite{BLT},
\cite{Ton}) and from orbit matrices (\cite{CrnEganS}, \cite{CRSR},
\cite{CrnRuk}, \cite{CrnDR}, \cite{HarTon}) have been extensively
studied. In \cite{0}, Tonchev described some extensions of incidence
matrices of a block design which define self-orthogonal binary codes
and in \cite{2} authors study similar extensions of orbit matrices and
submatrices of orbit matrices of $1$-design in order to construct
binary self-orthogonal codes from the simple group $\mathbf{He}$. In
\cite{CrnMost} authors construct self-dual codes over finite fields by
extending incidence matrices and orbit matrices of a block designs. In
\cite{Harada} proved that all $G$-invariant binary self-orthogonal
codes are subcodes of the dual code of the code spanned by the sets of
fixed points of involutions of G and they classify all
$M_{11}$-invariant binary self-dual codes. 

An incidence structure $\mathcal{D} =(\mathcal{P},\mathcal{B},\mathcal{I})$, with point set $\mathcal{P}$, block set $\mathcal{B}$ and incidence $\mathcal{I}$ is a $t$-$(v,k,\lambda)$ design, if $|\mathcal{P}|=v$, every block $B \in \mathcal{B}$ is incident with precisely $k$ points, and every $t$ distinct points are together incident with precisely $\lambda$ blocks. An incidence matrix of a $t$-$(v,k,\lambda)$ design $\mathcal{D}$ with $b$ blocks is a $b\times v$ matrix $M=[m_{i,j}],$ where $m_{i,j}=1$ if point $P_j$ is incident with block $B_i,$ and $0$ otherwise.
A $t$-$(v,k,\lambda)$ design is called {\bf weakly self-orthogonal} if all the block intersection numbers have the same parity. A design is {\bf self-orthogonal} if it is weakly self-orthogonal and if the block intersection numbers and the block size are even numbers. 

An isomorphism from one design to other is a bijective mapping of points to points and blocks to blocks which preserves incidence. An isomorphism from a design 
$\mathcal D$ onto itself is called an automorphism of $\mathcal{D}$. The set of all automorphisms of $\mathcal{D}$ forms its full automorphism group denoted by ${\rm Aut}(\mathcal D)$. 

The  {\bf code $C_{\mathbb{F}}(\mathcal{D})$ of the design} $\mathcal{D}$ over the finite field $\mathbb{F}$ is the space spanned by the incidence vectors of
the blocks over  $\mathbb{F}$. A code $C$ over field of order $q,$ length $n$, dimension $k$, and minimum weight $d$, is denoted by $[n,k,d]_q$ code. If $q=2,$ we denote code $C$ by $[n,k,d]$ code. The  {\bf dual} code $C^\perp$ is the orthogonal under the standard inner product, i.e.\ $C^{\perp} = \{ v \in \mathbb{F}^n \mid
v\cdot c=0 {\rm \ for\ all\ } c \in C \}$. 
A code $C$ is {\bf self-orthogonal} if $C \subseteq C^\perp$. If $\mathcal{D}$ is a self-orthogonal design then the binary code of the design $\mathcal{D}$ is self-orthogonal. The incidence matrix $M$ of a weakly self-orthogonal design can be extended in a certain way in order to span self-orthogonal code. The $i$-th row of the matrix $M$ will be denoted by $M[i].$
The all-one vector will be denoted by $\bf{1}$, and is the constant vector with all coordinate entries equal to 1. Two linear codes are {\bf equivalent} if they can be obtained from one another by multiplication of the coordinate positions by non-zero field elements or by permuting the coordinate positions. An {\bf automorphism} of a code $C$ is an isomorphism from $C$ to $C$. The full automorphism group will be denoted by ${\rm Aut}(C)$. If code $C_{\mathbb{F}}(\mathcal{D})$ is a linear code of a design $\mathcal D$ over a finite field $\mathbb{F},$ then the full automorphism group of $\mathcal{D}$ is contained in the full automorphism group of code $C_{\mathbb{F}}(\mathcal{D})$. Designs are obtained from transitive permutation representations of Mathieu group $M_{11}.$

$M_{11}$ is simple group of order $7920$ which has $39$ non-equivalent transitive permutation representations.
Among others, lattice of $M_{11}$ is consisted of $1$ subgroup of index $22$, $1$ subgroup of index $55$, $1$ subgroup of index $66$, $3$ subgroups of index $110,$ $2$subgroups of index $132$, $1$ subgroup of index $144$ and $1$ subgroup of index $165$. Subgroup of $M_{11}$ with largest index has index $3960.$
Using mentioned subgroups we obtained transitive permutation representations of $M_{11}$ on $22, 55, 66, 110, 132, 144$ and $165$ points.

In Section \ref{construction}, we generalize the construction (in \cite{0}) of binary self-orthogonal codes using incidence matrix of weakly self-orthogonal design to obtain self-orthogonal codes over an arbitrary finite field. In Section \ref{orbitne}, we extend the construction (in \cite{2}) of binary self-orthogonal codes using orbit matrices of self-orthogonal designs and weakly self-orthogonal designs such that block size is even and block intersection numbers are odd to all weakly self-orthogonal designs. Also, we generalize mentioned construction in order to obtain self-orthogonal codes over an arbitrary finite field. Finally, in Section \ref{results} we give tables of constructed binary codes obtained from weakly self-orthogonal designs invariant under an action of $M_{11},$ constructed using construction described in \cite{1}. Due to computational limitations, some codes are listed without their minimum distance and/or their automorphism groups. 

\section{Self-orthogonal codes obtained from WSO designs}\label{construction}
Let $\mathcal{D}$ be weakly self-orthogonal design and let $M$ be its $b\times v$ incidence matrix. Using suitable extension of $M$ one can obtain self-orthogonal binary code $C.$
\begin{thm}\label{tmTonchev}
Let $\mathcal{D}$ be weakly self-orthogonal design and let $M$ be its $b\times v$ incidence matrix.
\begin{enumerate}
	\item If $\mathcal{D}$ is a self-orthogonal design, than $C(\mathcal{D})$ is a binary self-orthogonal code.
	\item If $\mathcal{D}$ is such that $k$ is even and the block intersection numbers are odd, then the matrix $[I_b,M, \textbf{1}],$ generates a binary self-orthogonal code.
	\item If $\mathcal{D}$ is such that $k$ is odd and the block intersection numbers are even, then the matrix $[I_b,M],$ generates a binary self-orthogonal code.
	\item If $\mathcal{D}$ is such that $k$ is odd and the block intersection numbers are odd, then the matrix  $[M,\textbf{1}],$ generates a binary self-orthogonal code.
\end{enumerate}
\end{thm}
\textbf{Proof:} \cite{0}.
\begin{flushright}
$\blacksquare$ 
\end{flushright}
\noindent
Previous theorem can be generalized to obtain self-orthogonal codes over finite field $\mathbb{F}_q,$ where $q$ is prime power.

\begin{thm}\label{Tonchev_generalizacija}
Let $q=p^l$ be prime power and $\mathbb{F}_q$ a finite field of order $q.$ Let $\mathcal{D}$ be $1-(v,k,r)$ design and let $M$ be its $b\times v$ incidence matrix. Let $\mathcal{D}$ be such that $k \jmod{a}{p}$ and $|B_i\cap B_j| \jmod{d}{p},$ for all $i,j\in\{1,\ldots,b\}, \ i\neq j,$ where $B_i$ and $B_j$ are two blocks of the design $\mathcal{D}$.
\begin{enumerate}
	\item If $a=d=0,$ than $M$ generates a self-orthogonal code over $\mathbb{F}_q.$
	
	\item If $a=0$ and $d\neq 0$, then the matrix $[ \sqrt{d}\cdot I_b, M, \sqrt{-d}\cdot \textbf{1}]$ generates a $b$-dimensional self-orthogonal code over $\mathbb{F},$ where  $\mathbb{F}= \mathbb{F}_q$ if $d$ and $-d$ are quadratic residues modulo $q$, and $\mathbb{F}= \mathbb{F}_{q^2}$ otherwise.
	
	\item If $a\neq 0$ and $d= 0$, then the matrix $[ \sqrt{-a}\cdot I_b,M]$ generates a $b$-dimensional self-orthogonal code over $\mathbb{F},$ where $\mathbb{F}= \mathbb{F}_q$ if $-a$ is quadratic residue modulo $q$, and $\mathbb{F}= \mathbb{F}_{q^2}$ otherwise. If $b=v,$ obtained code is self-dual.
	
	\item If $a\neq 0$ and $d\neq 0$, there are two cases:
	\begin{enumerate}
				\item if $a=d,$ then the matrix $[ M, \sqrt{-a}\cdot \textbf{1}]$ generates a self-orthogonal code over $\mathbb{F},$ where  $\mathbb{F}= \mathbb{F}_q$ if $-a$ is quadratic residue modulo $q$, and $\mathbb{F}= \mathbb{F}_{q^2}$ otherwise, and
		
				\item if $a\neq d$, then the matrix $[ \sqrt{d-a}\cdot  I_b, M, \sqrt{-d}\cdot \textbf{1}]$ generates a $b$-dimensional self-orthogonal code over $\mathbb{F},$ where  $\mathbb{F}= \mathbb{F}_q$ if $d-a$ and $-d$ are quadratic residues modulo $q$, and $\mathbb{F}= \mathbb{F}_{q^2}$ otherwise.
	\end{enumerate} 
\end{enumerate}
\end{thm}
\textbf{Proof:}
\begin{enumerate}
		\item Since $a=0$, it follows that $$M[i]\cdot M[i]=k\jmod{a}{p}\jmod{0}{p},$$ and since $d=0$ it follows that $$M[i]\cdot M[j]=|B_i\cap B_j| \jmod{d}{p} \jmod{0}{p},$$ for all $i,j\in\{1,\ldots,b\}, \ i\neq j.$ 
		We conclude that $M$ generates a self-orthogonal code over $\mathbb{F}_q.$
		
		\item Let $A=[ \sqrt{d}\cdot I_b, M, \sqrt{-d}\cdot\textbf{1}].$ Then $$ A[i]\cdot A[i]=d+k-d \jmod{0}{p}$$ and $$ A[i]\cdot A[j]=|B_i\cap B_j|-d \jmod{d-d}{p} \jmod{0}{p},$$ for all $i,j\in\{1,\ldots,b\}, \ i\neq j.$ We conclude that $A$ generates a self-orthogonal code over $\mathbb{F},$ where  $\mathbb{F}= \mathbb{F}_q$ if $-d$ is quadratic residue modulo $q$, and $\mathbb{F}= \mathbb{F}_{q^2}$ otherwise. \footnote{Elements which are equal to 0 in $\mathbb{F}_q$ are also equal to 0 in $\mathbb{F}_{q^2},$ since elements of $\mathbb{F}_q$ are polynomials in $\mathbb{F}_{q^2}$ of degree at most 1 with coefficients in $\mathbb{F}_q$.}
		
		\item Let $A=[ \sqrt{-a}\cdot I_b,M].$ Then $$ A[i]\cdot A[i]=-a+k \jmod{0}{p}$$ and $$ A[i]\cdot A[j]=|B_i\cap B_j| \jmod{0}{p},$$ for all $i,j\in\{1,\ldots,b\}, \ i\neq j.$ We conclude that $A$ generates a $b$-dimensional self-orthogonal code over $\mathbb{F},$ where  $\mathbb{F}= \mathbb{F}_q$ if $-d$ is quadratic residue modulo $q$, and $\mathbb{F}= \mathbb{F}_{q^2}$ otherwise. If $b=v,$ dimension of a code is half of its length, so the code is self-dual.
  	
		\item \begin{enumerate}
				\item Let $A=[ M, \sqrt{-a}\cdot \textbf{1}].$ Then $$ A[i]\cdot A[i]=k-a \jmod{0}{p}$$ and $$ A[i]\cdot A[j]=|B_i\cap B_j|-a \jmod{0}{p},$$ for all $i,j\in\{1,\ldots,b\}, \ i\neq j.$ We conclude that $A$ generates a self-orthogonal code over $\mathbb{F},$ where  $\mathbb{F}= \mathbb{F}_q$ if $-d$ is quadratic residue modulo $q$, and $\mathbb{F}= \mathbb{F}_{q^2}$ otherwise.
		
				\item Let $A=[ \sqrt{d-a}\cdot I_b, M, \sqrt{-d}\cdot\textbf{1}].$ Then $$A[i]\cdot A[i]=d-a+k-d \jmod{0}{p}$$ and $$ A[i]\cdot A[j]=|B_i\cap B_j|-d \jmod {0}{p},$$ for all $i,j\in\{1,\ldots,b\}, \ i\neq j.$ We conclude that $A$ generates a $b$-dimensional self-orthogonal code over $\mathbb{F},$ where $\mathbb{F}= \mathbb{F}_q$ if $-d$ is quadratic residue modulo $q$, and $\mathbb{F}= \mathbb{F}_{q^2}$ otherwise.
	\end{enumerate} 
\end{enumerate}
\begin{flushright}
$\blacksquare$
\end{flushright}

\section{Codes from orbit matrices of weakly self-orthogonal designs}\label{orbitne}

\begin{rmk} \label{remark}
Let $\mathcal{D}$ be $1-(v,k,r)$ design and $G$ be an automorphism group of the design. Let $v_1=|\mathcal{V}_1|,\ldots,v_n=|\mathcal{V}_n|$ be the sizes of point orbits and $b_1=|\mathcal{B}_1|,\ldots,b_m=|\mathcal{B}_m|$ be the sizes of block orbits under the action of the group $G.$ We define an orbit matrix under the action of $G$ as $m\times n$ matrix:
$$\begin{bmatrix}
a_{1,1} & a_{1,2} & \ldots & a_{1,n} \\
a_{2,1} & a_{2,2} & \ldots & a_{2,n} \\
\vdots & \vdots & \ddots & \vdots \\
a_{m,1} & a_{m,2} & \ldots & a_{m,n} \\
\end{bmatrix},$$
where $a_{i,j}$ is the number of points of the orbit $\mathcal{V}_j$ incident with a block of the orbit $\mathcal{B}_i.$ 
\\
It is easy to see that the orbit matrix is well defined and that $k=\sum_{j=1}^n a_{i,j}.$
\\
For $x\in\mathcal{B}_s$, by counting the incidence pairs $(P,x')$ such that $x'\in \mathcal{B}_t$ and $P$ is incident with the block $x,$ we obtain $\sum_{x'\in\mathcal{B}_t} |x\cap x'|=\sum_{j=1}^n \frac{b_t}{v_j}a_{s,j}a_{t,j}.$
\end{rmk}

\begin{rmk} \label{remarkmodp}
Let $\mathcal{D}$ be $1-(v,k,r)$ design such that $k \jmod{a}{p}$ and $|B_i\cap B_j| \jmod{d}{p},$ for all $i,j\in\{1,\ldots,b\}, \ i\neq j,$ where $B_i$ and $B_j$ are two blocks of the design $\mathcal{D}.$ Let $G$ be an automorphism group of the design which acts on $\mathcal{D}$ with $n$ point orbits of length $w$ and block orbits of length $b_1,b_2,\ldots,b_m,$ and let $O$ be an orbit matrix of a design $\mathcal{D}$ under the action of the group $G$. 
From previous remark, for $x\in\mathcal{B}_s$ and $s\neq t$ one concludes that
\begin{align*}
\frac{b_t}{w} O[s]\cdot O[t]&=\sum_{j=1}^n \frac{b_t}{w}a_{s,j}a_{t,j}\\
														&=\sum_{x'\in\mathcal{B}_t} |x\cap x'|.
\end{align*}
It follows that $$\frac{b_t}{w}O[s]\cdot O[t] \jmod{b_t d}{p}.$$

Similar, for $x\in\mathcal{B}_s$ one can conclude that
\begin{align*}
\frac{b_t}{w} O[s]\cdot O[s]&=\sum_{x'\in\mathcal{B}_s} |x\cap x'|\\
&=|x\cap x|+\sum_{x'\in\mathcal{B}_s,\ x\neq x'} |x\cap x'|.																	
\end{align*}
It follows that $$\frac{b_s}{w} O[s]\cdot O[s] \jmod{a+(b_s-1)d}{p}.$$
\end{rmk}

Let $\mathcal{D}$ be $1-(v,k,r)$ design and $G$ an automorphism group of $\mathcal{D}$ acting with $f_1$ fixed points and $n$ point orbits of length $q>1,$ and with $f_2$ fixed blocks and $m$ block orbits of length $q.$ We define matrices $OM1$ and $OM2$ to be, respectively, the matrices 
$$\begin{bmatrix}
 a_{1,1} & a_{1,2} & \ldots & a_{1,f_1} \\
a_{2,1} & a_{2,2} & \ldots & a_{2,f_1}\\
\vdots & \vdots & \ddots & \vdots \\
 a_{f_2,1} & a_{f_2,2} & \ldots & a_{f_2,f_1} \\
\end{bmatrix}$$
 and 
$$\begin{bmatrix}
 a_{f_2+1,f_1+1} & a_{f_2+1,f_1+2} & \ldots & a_{f_2+1,f_1+n} \\
a_{f_2+2,f_1+1} & a_{f_2+2,f_1+2} & \ldots & a_{f_2+2,f_1+n}\\
 \vdots & \vdots & \ddots & \vdots \\
 a_{f_2+m,f_1+1} & a_{f_2+m,f_1+2} & \ldots & a_{f_2+m,f_1+n}\\
\end{bmatrix},$$
where the columns $1,2,\ldots, f_1$ correspond to the fixed points and the rows $1,2,\ldots,f_2$ correspond to the fixed blocks.

\begin{rmk}\label{remark_sum}
\begin{itemize}
	\item[a)] If $B_1$ and $B_2$ are blocks fixed under the action of group $G$ on the design, then $B_1, \ B_2$ and $B_1\cap B_2$ are unions of some $G$-orbits of the point set.
	\item[b)] Let $\mathcal{B}_t$ and $\mathcal{B}_s$ be block orbits of size $q$ under the action of the group $G$ on the design. It follows from Remark \ref{remark} that $$\sum_{x'\in\mathcal{B}_t} |x\cap x'|=\sum_{j=1}^{f_1} \frac{b_t}{v_j}a_{s,j}a_{t,j}+\sum_{j=f_1+1}^{f_1+n} \frac{b_t}{v_j}a_{s,j}a_{t,j}=q\sum_{j=1}^{f_1}  a_{s,j}a_{t,j}+\sum_{j=f_1+1}^{f_1+n}a_{s,j}a_{t,j}.$$
\end{itemize}
\end{rmk}


\subsection{Codes from orbit matrices of self-orthogonal 1-designs}
\begin{thm}\label{tmsl1}
Let $\mathcal{D}$ be a self-orthogonal 1-design and $G$ be an automorphism group of the design which acts on $\mathcal{D}$ with $n$ point orbits of length $w$ and block orbits of length $b_1,b_2,\ldots,b_m$ such that $b_i=2^o\cdot b_i', \ w=2^u\cdot w',\ o\leq u, 2\notdivides b_i'$ and $2\notdivides w'$, for $i\in\{1,\ldots,m\}.$
Then the binary code spanned by the rows of orbit matrix of the design $\mathcal{D}$ (under the action of the group $G$) is a self-orthogonal code of length $\frac{v}{w}$.
\end{thm}
\textbf{Proof:}
In \cite{2}. 
\begin{flushright}
$\blacksquare$
\end{flushright}
Using the following theorem, one can use orbit matrices of a self-orthogonal design to construct self-orthogonal codes over $\mathbb{F}_q,$ where $q$ is prime power.

\begin{thm}
Let $q=p^l$ be prime power and $\mathbb{F}_q$ a finite field of order $q.$ Let $\mathcal{D}$ be a $1-(v,k,r)$ design such that $k \jmod{0}{p}$ and $|B_i\cap B_j| \jmod{0}{p},$ for all $i,j\in\{1,\ldots,b\}, \ i\neq j,$ where $B_i$ and $B_j$ are two blocks of the design $\mathcal{D},$ and let $G$ be an automorphism group of the design which acts on $\mathcal{D}$ with $n$ point orbits of length $w$ and $m$ block orbits of length $w.$
Then the binary linear code spanned by the rows of orbit matrix of the design $\mathcal{D}$ (under the action of the group $G$) is a self-orthogonal code over $\mathbb{F}_q$ of length $\frac{v}{w}$ with dimension equal to $\rm{rank} O.$
\end{thm}
\textbf{Proof:}
Let $O$ be the orbit matrix of a design $\mathcal{D}$ under the action of $G$. Form Remark \ref{remarkmodp} it follows that $O[s]\cdot O[t] \jmod{0}{p},$ and $O[s]\cdot O[s]\jmod{0}{p},$ for all $s,t\in\{1,\ldots,m\}.$
\begin{flushright}
$\blacksquare$
\end{flushright}
%
%
%
%
%
%
%
\begin{thm}\label{fixpoint1}
Let $\mathcal{D}$ be a self-orthogonal $1-(v,k,r)$ design and $G$ be an automorphism group of $\mathcal{D}$ acting with $f_1$ fixed points and $n$ point orbits of length $2,$ and with $f_2$ fixed blocks and $m$ block orbits of length $2.$ Then
\begin{enumerate}
	\item the binary linear code spanned by the matrix $OM1$ is a binary self-orthogonal code of length $f_1$ and dimension $\rm{rank}(OM1);$
	
	\item the binary linear code spanned by the matrix $OM2$ is a binary self-orthogonal code of length $n$ and dimension $\rm{rank}(OM2).$
\end{enumerate}
\end{thm}
\textbf{Proof:} In \cite{2}.
\begin{flushright}
$\blacksquare$
\end{flushright}
%
%
%
%
We generalize previous theorem to obtain self-orthogonal codes over $\mathbb{F}_q,$ where $q$ is a prime power.
\begin{thm}\label{ftsl1_q}
Let $q=p^l$ be prime power and $\mathbb{F}_q$ a finite field of order $q.$ Let $\mathcal{D}$ be a $1-(v,k,r)$ design such that $k \jmod{0}{p}$ and $|B_i\cap B_j| \jmod{0}{p},$ for all $i,j\in\{1,\ldots,b\}, \ i\neq j,$ where $B_i$ and $B_j$ are two blocks of the design $\mathcal{D},$ and let $G$ be an automorphism group of the design which acts on $\mathcal{D}$ with $f_1$ fixed points and $n$ point orbits of length $p^\alpha,\ 1\leq \alpha\leq l$, and with $f_2$ fixed blocks and $m$ block orbits of length $p^\alpha$. 
Then 
 \begin{enumerate}
	\item the linear code spanned by the matrix $OM1$ is a self-orthogonal $[f_1,\rm{rank}(OM1)]$ code over the field $\mathbb{F}_q.$ 
		\item the linear code spanned by the matrix $OM2$ is a self-orthogonal $[n,\rm{rank}(OM2)]$ code over the field $\mathbb{F}_q.$ 
\end{enumerate}
\end{thm}
\textbf{Proof:}
 \begin{enumerate}
	\item Since $k\jmod{0}{p}$, each block contains $p\cdot \beta$ fixed points and since $|B_i\cap B_j| \jmod{0}{p},$ for all $i,j\in\{1,\ldots,b\}, \ i\neq j,$ intersection of every two blocks contains $q\cdot \gamma$ fixed points. Therefore, it follows that the matrix $OM1$ spans a self-orthogonal code over the field $\mathbb{F}_q.$

	\item From Remark \ref{remark_sum}, for $s\neq t,$ since $\mathcal{B}_t$ is block orbit of size $p^\alpha$, it follows that
	$$\sum_{x'\in\mathcal{B}_t} |x\cap x'| \jmod{0}{p},$$ 
	and for $s=t$ it follows that
	$$\sum_{x'\in\mathcal{B}_s} |x\cap x'|=|x\cap x|+ \sum_{x'\in\mathcal{B}_s, \\ x\neq x'}|x\cap x'| \jmod{0+(p^\alpha-1)0}{p}.$$
	We conclude that 
	$OM2[s]\cdot OM2[t] \jmod{0}{p},$ for $s\neq t$ and $OM2[s]\cdot OM2[s] \jmod{0}{p}$ and that the binary code spanned by the matrix $OM2$ is self-orthogonal over the field $\mathbb{F}_q.$
\end{enumerate}
\begin{flushright}
$\blacksquare$
\end{flushright}


\subsection{Codes from orbit matrices of extended weakly self-orthogonal 1-designs with $k$ even and odd block intersection numbers}
\begin{thm}\label{tmsl2}
Let $\mathcal{D}$ be a weakly self-orthogonal 1-design such that $k$ is even and the block intersection numbers are odd, and let $G$ be an automorphism group of the design which acts on $\mathcal{D}$ with $n$ point orbits of length $w$ and block orbits of length $b_1,b_2,\ldots,b_m$ such that $b_i=2^o\cdot b_i', \ w=2^u\cdot w',\ o\leq u, 2\notdivides b_i'$ and $2\notdivides w'$ for $i\in\{1,\ldots,m\}.$ Let $O$ be the orbit matrix of $\mathcal{D}$ under action of the group $G.$ 

\begin{itemize}
		\item[a)] If $o=u=0,$ then the binary linear code spanned by the rows of the matrix $[I_m, O, \textbf{1}]$ is a self-orthogonal code of the length $m+\frac{v}{w}+1$ and dimension $m.$ 
		\item[b)] If $o\geq 1$ and $o=u$ then the binary linear code spanned by the rows of the matrix $[I_m, O]$ is a self-orthogonal code of the length $m+\frac{v}{w}$ and dimension $m.$ If $m=n,$ the obtained code is self-dual.
		\item[c)] If $o<u,$ then binary linear code spanned by the rows of the matrix $O$ is a self-orthogonal code of the length $\frac{v}{w}$ and dimension $\rm{rank}O.$
\end{itemize}
\end{thm}
\textbf{Proof:}
\begin{itemize}
		\item[a)] Since $w, b_1, \ldots, b_m$ are all odd numbers, it follows from Remark \ref{remarkmodp} that 
		$$O[s]\cdot O[t] \jmod{1}{2},\ \ O[s]\cdot O[s] \jmod{0}{2},$$ for all $s,t\in\{1,\ldots,m\},\ s\neq t.$ 
		\item[b)] Since $w, b_1, \ldots, b_m$ are all even numbers and $\frac{b_t}{w}=\frac{2^o\cdot b_t'}{2^o\cdot w'}=\frac{b_t'}{w'},$ for every $t\in\{1,\ldots,m\},$ it follows from Remark \ref{remarkmodp} that $$O[s]\cdot O[t] \jmod{0}{2},\ \ O[s]\cdot O[s] \jmod{1}{2},$$ for all $s,t\in\{1,\ldots,m\}.$ 
		\item[c)] Since $o<u$ and since $\displaystyle\frac{b_t}{w}O[s]\cdot O[t]=\frac{b_t'}{2^{u-o}w'}O[s]\cdot O[t]$ is a positive integer from Remark \ref{remarkmodp}, it follows that $2\mid O[s]\cdot O[t],$ for all $s,t\in\{1,\ldots,m\}.$
\end{itemize}
\begin{flushright}
$\blacksquare$
\end{flushright}
\noindent
One can generalize previous theorem to construct a self-orthogonal code over $\mathbb{F}_q,$ where $q$ is a prime power. The following theorem describes the construction.

\begin{thm}\label{thmsl2_q}
Let $q=p^l$ be prime power and $\mathbb{F}_q$ a finite field of order $q.$ Let $\mathcal{D}$ be $1-(v,k,r)$ design such that $k \jmod{0}{p}$ and $|B_i\cap B_j| \jmod{d}{p},$ for all $i,j\in\{1,\ldots,b\}, \ i\neq j,$ where $B_i$ and $B_j$ are two blocks of the design $\mathcal{D},$ and let $G$ be an automorphism group of the design which acts on $\mathcal{D}$ with $n$ point orbits of length $w$ and $m$ block orbits of length $w$ and let $O$ be the orbit matrix of $\mathcal{D}$ under action of the group $G.$ 

\begin{itemize}
		\item[a)] If $p\mid w$, then linear code spanned by the rows of the matrix $A=[\sqrt{d} I_m, O]$ is a self-orthogonal $[m+n,m]$ code over the field $\mathbb{F},$ where $\mathbb{F}=\mathbb{F}_q$ if $d$ is square root in $\mathbb{F}_q,$ and $\mathbb{F}=\mathbb{F}_{q^2}$ otherwise. If $m=n$, tha obtained code is self-dual.
		\item[b)] If $p\mid w-1$, then linear code spanned by the rows of the matrix $A=[\sqrt{wd} I_m, O, \sqrt{-wd}\textbf{1}]$ is a self-orthogonal $[m+n+1,m]$ code over the field $\mathbb{F},$ where $\mathbb{F}=\mathbb{F}_q$ if $wd$ and $-wd$ are square roots in $\mathbb{F}_q,$ and $\mathbb{F}=\mathbb{F}_{q^2}$ otherwise.
		\item[c)] If $p\notdivides  w$ and $p\notdivides  w-1$, then linear code spanned by the rows of the matrix $A=[\sqrt{d} I_m, O, \sqrt{-wd}\textbf{1}]$ is a self-orthogonal $[m+n+1,m]$ code over the field $\mathbb{F},$ where $\mathbb{F}=\mathbb{F}_q$ if $d$ and $-wd$ are square roots in $\mathbb{F}_q,$ and $\mathbb{F}=\mathbb{F}_{q^2}$ otherwise.
\end{itemize}
\end{thm}
\textbf{Proof:}
It follows from Remark \ref{remarkmodp} that $O[s]\cdot O[t] \jmod{wd}{p},$ for $s\neq t$ and $ O[s]\cdot O[s] \jmod{(w-1)d}{p}.$
\begin{itemize}
		\item[a)] If $p\mid w,$ then $A[s]\cdot A[t]=O[s]\cdot O[t] \jmod{0}{p}$ and $A[s]\cdot A[s]=-d+O[s]\cdot O[s] \jmod{0}{p}.$
		\item[b)] If $p\mid w-1,$ then $A[s]\cdot A[t]=O[s]\cdot O[t]-wd \jmod{0}{p}$ and $A[s]\cdot A[s]=wd+O[s]\cdot O[s]-wd \jmod{0}{p}.$
		\item[c)] If $p\notdivides  w$ and $p\notdivides  w-1,$ then $A[s]\cdot A[t]=O[s]\cdot O[t]-wd \jmod{0}{p}$ and $A[s]\cdot A[s]=d+O[s]\cdot O[s]-wd \jmod{0}{p}.$
\end{itemize}
\begin{flushright}
$\blacksquare$
\end{flushright}




\begin{thm}\label{ft1}
Let $\mathcal{D}$ be a weakly self-orthogonal $1-(v,k,r)$ design such that $k$ is even and block intersection numbers are odd. Let $G$ be an automorphism group of $\mathcal{D}$ acting with $f_1$ fixed points and $n$ point orbits of length $2,$ and with $f_2$ fixed blocks and $m$ block orbits of length $2.$ 
Then
\begin{enumerate}
	\item the binary linear code spanned by the matrix $[I_{f_2}, OM1, \textbf{1}]$ is a self-orthogonal code of length $f_2+f_1+1$ and dimension $f_2$,
	\item the binary linear code spanned by the matrix $[I_m, OM2]$ is a self-orthogonal code of length $m+n$ and dimension $m$. If $m=n,$ the obtained code is self-dual.
\end{enumerate}
\end{thm}
\textbf{Proof:}
 \begin{enumerate}
	\item Since $k$ is even, each block contains even number of fixed points and since block intersection numbers are odd, intersection of every two blocks contains odd number of fixed points. Therefore, it follows that the matrix 
$[I_{f_2}, OM1, \textbf{1}]$ spans a binary self-orthogonal code.

	\item For $s\neq t,$ since $\mathcal{B}_t$ is block orbit of size 2, it follows form Remark \ref{remark_sum} that 
	$$\sum_{x'\in\mathcal{B}_t} |x\cap x'| \jmod{2\cdot 1}{2} \jmod{0}{2},$$
	and for $s=t,$
	$$\sum_{x'\in\mathcal{B}_s} |x\cap x'|=|x\cap x|+ \sum_{x'\in\mathcal{B}_s, \\ x\neq x'}|x\cap x'| \jmod{0+1}{2}.$$
	We conclude that 
	$\sum_{j=f_1+1}^{f_1+n}a_{s,j}a_{t,j} \jmod{0}{2},$ for $s\neq t$ and $\sum_{j=f_1+1}^{f_1+n}a_{s,j}a_{s,j} \jmod{1}{2}$ and that the binary code spanned by the matrix $[I_m, OM2]$ is self-orthogonal. If $m=n,$ dimension of the code is half of its length, so the code is self-dual.
\end{enumerate}
\begin{flushright}
$\blacksquare$
\end{flushright}
We generalize previous theorem to obtain self-orthogonal codes over $\mathbb{F}_q,$ where $q$ is a prime power.

\begin{thm}\label{ftsl2_q}
Let $q=p^l$ be prime power and $\mathbb{F}_q$ a finite field of order $q.$ Let $\mathcal{D}$ be $1-(v,k,r)$ design such that $k \jmod{0}{p}$ and $|B_i\cap B_j| \jmod{d}{p},$ for all $i,j\in\{1,\ldots,b\}, \ i\neq j,$ where $B_i$ and $B_j$ are two blocks of the design $\mathcal{D},$ and let $G$ be an automorphism group of the design which acts on $\mathcal{D}$ with $f_1$ fixed points and $n$ point orbits of length $p^\alpha, \ 1\leq\alpha\leq l$, and with $f_2$ fixed blocks and $m$ block orbits of length $p^\alpha$. 
Then 
 \begin{enumerate}
	\item linear code spanned by the matrix $[\sqrt{d}\cdot I_{f_2}, OM1, \sqrt{-d}\cdot\textbf{1}]$ is a self-orthogonal $[f_2+f_1+1,f_2]$ code over the field $\mathbb{F},$ where $\mathbb{F}=\mathbb{F}_q$ if $d$ and $-d$ are square roots in $\mathbb{F}_q,$ and $\mathbb{F}=\mathbb{F}_{q^2}$ otherwise.
	\item linear code spanned by the matrix $[\sqrt{d}\cdot I_m, OM2]$ is a self-orthogonal $[m+n,m]$ code over the field $\mathbb{F},$ where $\mathbb{F}=\mathbb{F}_q$ if $d$ is square root in $\mathbb{F}_q,$ and $\mathbb{F}=\mathbb{F}_{q^2}$ otherwise. If $m=n,$ the obtained code is self-dual.
\end{enumerate}
\end{thm}
\textbf{Proof:}
 \begin{enumerate}
	\item Since $k\jmod{0}{p}$, each block contains $p\cdot \beta$ fixed points and since $|B_i\cap B_j| \jmod{d}{p},$ for all $i,j\in\{1,\ldots,b\}, \ i\neq j,$, intersection of every two blocks contains $p\cdot \gamma+d$ fixed points. Therefore, it follows that the matrix $[\sqrt{d}\cdot I_m, OM1, \sqrt{-d}\cdot\textbf{1}]$ spans a self-orthogonal code of length $f_2+f_1+1$ over the field $\mathbb{F}.$

	\item For $s\neq t,$ since $\mathcal{B}_t$ is block orbit of size $p^\alpha$, it follows from Remark \ref{remark_sum} that 
	$$\sum_{x'\in\mathcal{B}_t} |x\cap x'| \jmod{p^\alpha d}{p} \jmod{0}{p},$$
	and if $s=t$ it follows that
	$$\sum_{x'\in\mathcal{B}_s} |x\cap x'|=|x\cap x|+ \sum_{x'\in\mathcal{B}_s, \\ x\neq x'}|x\cap x'| \jmod{0+(p^\alpha-1)d}{p}\jmod{-d}{p}.$$
	We conclude that 
	$OM2[s]\cdot OM2[t] \jmod{0}{p},$ for $s\neq t$ and $OM2[s]\cdot OM2[s] \jmod{(p^\alpha-1)d}{p}$ and that the binary code spanned by the matrix $[\sqrt{d}\cdot I_m, OM2]$ is self-orthogonal over the field $\mathbb{F}.$
\end{enumerate}
\begin{flushright}
$\blacksquare$
\end{flushright}


\subsection{Codes from orbit matrices of extended weakly self-orthogonal 1-designs with $k$ odd and even block intersection numbers}
\begin{thm}\label{tmsl3}
Let $\mathcal{D}$ be a weakly self-orthogonal 1-design such that $k$ is odd and the block intersection numbers are even and $G$ be an automorphism group of the design which acts on $\mathcal{D}$ with $n$ point orbits of length $w$ and block orbits $b_1,b_2,\ldots,b_m$ such that $b_i=2^o\cdot b_i', \ w=2^u\cdot w',\ o\leq u, 2\notdivides b_i'$ and $2\notdivides w'$, for $i\in\{1,\ldots,m\}.$ Let $O$ be the orbit matrix of $\mathcal{D}$ under action of the group $G.$ 
\begin{itemize}
		\item[a)] If $o=u$, then he binary linear code spanned by the rows of matrix $[I_m, O]$ is a self-orthogonal code of length $m+\frac{v}{w}$ and dimension $m.$ If $m=n,$ the obtained code is self-dual.
		\item[b)] If $o<u$, then he binary linear code spanned by the rows of matrix $O$ is a self-orthogonal code of length $\frac{v}{w}$ and dimension $\rm{rank}O.$
\end{itemize}
\end{thm}
\textbf{Proof:} 
In \cite{2}. 
\begin{flushright}
$\blacksquare$
\end{flushright}
One can generalize previous theorem to construct a self-orthogonal code over $\mathbb{F}_q,$ where $q$ is a prime power. The following theorem describes the construction.
\begin{thm}
Let $q=p^l$ be prime power and $\mathbb{F}_q$ a finite field of order $q.$ Let $\mathcal{D}$ be $1-(v,k,r)$ design such that $k \jmod{a}{p}$ and $|B_i\cap B_j| \jmod{0}{p},$ for all $i,j\in\{1,\ldots,b\}, \ i\neq j,$ where $B_i$ and $B_j$ are two blocks of the design $\mathcal{D},$ and let $G$ be an automorphism group of the design which acts on $\mathcal{D}$ with $n$ point orbits of length $w$ and $m$ block orbits of length $w.$ 
Then the linear code spanned by the rows of matrix $A=[\sqrt{-a}I_m, O],$ where $O$ is the orbit matrix of the design $\mathcal{D}$ (under the action of the group $G$), is a self-orthogonal $[m+n,m]$ code over $\mathbb{F},$ where $\mathbb{F}=\mathbb{F}_q$ if $a$ is a square root modulo $q$, and $\mathbb{F}=\mathbb{F}_{q^2}$ otherwise. If $m=n,$ the obtained code is self-dual.
\end{thm}
\textbf{Proof:}
From Remark \ref{remarkmodp} it follows that $O[s]\cdot O[t]\jmod{0}{p},$ and $O[s]\cdot O[s]\jmod{a}{p},$ so the rows of $A=[\sqrt{-a}I_m, O]$ generates self-orthogonal code. If $m=n,$ the dimension of the code is half of its length, so the code is self-dual.  
\begin{flushright}
$\blacksquare$
\end{flushright}

\begin{thm}\label{ftsl3}
Let $\mathcal{D}$ be a weakly self-orthogonal $1-(v,k,r)$ design such that $k$ is odd and block intersection numbers are even. Let $G$ be an automorphism group of $\mathcal{D}$ acting with $f_1$ fixed points and $n$ point orbits of length $2,$ and with $f_2$ fixed blocks and $m$ block orbits of length $2.$ Then
\begin{enumerate}
	\item the binary linear code spanned by the matrix $[I_{f_2},OM1]$ is a binary self-orthogonal code of length $f_2+f_1$ and dimension $f_2;$
	\item the binary linear code spanned by the matrix $[I_{m}, OM2]$ is a self-orthogonal code of length $m+n$ and dimension $m.$
\end{enumerate}
\end{thm}
\textbf{Proof:} In \cite{2}. 
\begin{flushright}
$\blacksquare$
\end{flushright}
\noindent
We generalize previous theorem to obtain self-orthogonal codes over $\mathbb{F}_q,$ where $q$ is a prime power.

\begin{thm}\label{ftsl3_q}
Let $q=p^l$ be prime power and $\mathbb{F}_q$ a finite field of order $q.$ Let $\mathcal{D}$ be $1-(v,k,r)$ design such that $k \jmod{a}{p}$ and $|B_i\cap B_j| \jmod{0}{p},$ for all $i,j\in\{1,\ldots,b\}, \ i\neq j,$ where $B_i$ and $B_j$ are two blocks of the design $\mathcal{D},$ and let $G$ be an automorphism group of the design which acts on $\mathcal{D}$ with $f_1$ fixed points and $n$ point orbits of length $p^\alpha, \ 1\leq\alpha\leq l$, and with $f_2$ fixed blocks and $m$ block orbits of length $p^\alpha$. 
Then 
 \begin{enumerate}
	\item the linear code spanned by the matrix $[\sqrt{-a}\cdot I_{f_2}, OM1]$ is a self-orthogonal $[f_2+f_1, f_2]$ code over the field $\mathbb{F},$ where $\mathbb{F}=\mathbb{F}_q$ if $-a$ is square root in $\mathbb{F}_q,$ and $\mathbb{F}=\mathbb{F}_{q^2}$ otherwise.
	\item the linear code spanned by the matrix $[\sqrt{-a}\cdot I_m, OM2]$ is a self-orthogonal $[m+n,m]$ code over the field $\mathbb{F},$ where $\mathbb{F}=\mathbb{F}_q$ if $-a$ is square root in $\mathbb{F}_q,$ and $\mathbb{F}=\mathbb{F}_{q^2}$ otherwise. If $m=n,$ the obtained code is self-dual.
\end{enumerate}
\end{thm}
\textbf{Proof:}
 \begin{enumerate}
	\item Since $k\jmod{a}{p},$ each block contains $p\cdot \beta+a$ fixed points and since $|B_i\cap B_j| \jmod{0}{p},$ for all $i,j\in\{1,\ldots,b\}, \ i\neq j,$ intersection of every two blocks contains $q\cdot \gamma$ fixed points. Therefore, it follows that the matrix $[\sqrt{-a}\cdot I_m, OM1]$ spans a self-orthogonal code of length $f_2+f_1+1$ over the field $\mathbb{F}.$

	\item	For $s\neq t,$ since $\mathcal{B}_t$ is block orbit of size $p^\alpha$, it follows from Remark \ref{remark_sum} that 
	$$\sum_{x'\in\mathcal{B}_t} |x\cap x'|  \jmod{0}{p},$$
	and for $s=t$ it follows that
	$$\sum_{x'\in\mathcal{B}_s} |x\cap x'|=|x\cap x|+ \sum_{x'\in\mathcal{B}_s, \\ x\neq x'}|x\cap x'| \jmod{a+0}{p}.$$
	We conclude that 
	$OM2[s]\cdot OM2[t] \jmod{0}{p},$ for $s\neq t$ and $OM2[s]\cdot OM2[s] \jmod{a}{p}$ and that the binary code spanned by the matrix $[\sqrt{-a}\cdot I_m, OM2]$ is self-orthogonal over the field $\mathbb{F}.$ 
\end{enumerate}
\begin{flushright}
$\blacksquare$
\end{flushright}




\subsection{Codes from orbit matrices of extended weakly self-orthogonal 1-designs with $k$ odd and odd block intersection numbers}
\begin{thm}\label{tmsl4}
Let $\mathcal{D}$ be a weakly self-orthogonal 1-design such that $k$ is odd and the block intersection numbers are odd and $G$ be an automorphism group of the design which acts on $\mathcal{D}$ with $n$ point orbits of length $w$ and block orbits of length $b_1,b_2,\ldots,b_m$ such that $b_i=2^o\cdot b_i', \ w=2^u\cdot w',\ o\leq u, 2\notdivides b_i'$ and $2\notdivides w'$, for $i\in\{1,\ldots,m\}.$ Let $O$ be the orbit matrix of $\mathcal{D}$ under action of the group $G.$ 

\begin{itemize}
		\item[a)] If $o=u=0,$ then the binary linear code spanned by the rows of the matrix $[O, \textbf{1}]$ is a self-orthogonal code of the length $\frac{v}{w}+1$ and dimension $\rm{rank}O.$
		\item[b)] Otherwise, the binary linear code spanned by the rows of the matrix $O$ is a self-orthogonal code of the length $\frac{v}{w}$ and dimension $\rm{rank}O.$
\end{itemize}

\end{thm}
\textbf{Proof:} 
		\begin{itemize}
				\item[a)] If $o=u=0,$ it follows from Remark \ref{remarkmodp} that $$\frac{b_t}{w} O[s]\cdot O[t] \jmod{1}{2}$$ and $$\frac{b_s}{w} O[s]\cdot O[s] \jmod{1}{2}$$ for all $s,t\in\{1,\ldots,m\}.$ Since $w, b_1,\ldots,b_s$ are odd numbers, it follows that $O[s]\cdot O[t] \jmod{1}{2}$ and $O[s]\cdot O[s] \jmod{1}{2}.$
				\item[b)] If $o=u>1,$ it follows Remark \ref{remarkmodp} that $$\frac{b_t'}{w'} O[s]\cdot O[t] \jmod{0}{2}$$ for all $s,t\in\{1,\ldots,m\}.$ Since $w', b_1',\ldots,b_s'$ are odd numbers, it follows that $O[s]\cdot O[t] \jmod{0}{2}$ and $O[s]\cdot O[s] \jmod{0}{2}.$\\
									If $1<o<u,$ it follows Remark \ref{remarkmodp} that $$\frac{b_s'}{2^{u-o}\cdot w'} O[s]\cdot O[t] \jmod{0}{2}$$ for all $s,t\in\{1,\ldots,m\}$ and since $\displaystyle \frac{b_s'}{2^{u-o}\cdot w'} O[s]\cdot O[t]$ has to be a positive integer, it follows that $ O[s]\cdot O[t] \jmod{0}{2}$ and $O[s]\cdot O[s] \jmod{0}{2}.$
		\end{itemize}
\begin{flushright}
$\blacksquare$
\end{flushright}
One can generalize previous theorem to construct a self-orthogonal code over $\mathbb{F}_q,$ where $q$ is a prime power. The following theorem describes the construction.

\begin{thm}\label{thmsl4_q}
Let $q=p^l$ be prime power and $\mathbb{F}_q$ a finite field of order $q.$ Let $\mathcal{D}$ be $1-(v,k,r)$ design such that $k \jmod{a}{p}$ and $|B_i\cap B_j| \jmod{d}{p},$ for all $i,j\in\{1,\ldots,b\}, \ i\neq j,$ where $B_i$ and $B_j$ are two blocks of the design $\mathcal{D},$ and let $G$ be an automorphism group of the design which acts on $\mathcal{D}$ with $n$ point orbits of length $w$ and $m$ block orbits of length $w$ and let $O$ be the orbit matrix of $\mathcal{D}$ under action of the group $G.$ 

\begin{itemize}

	\item If $a=d$ we differ two cases.
	
	\begin{itemize}
			\item[a)] If $p\mid w$, then linear code spanned by the rows of the matrix $O$ is a self-orthogonal $[m,\rm{rank}O]$ code over the field $\mathbb{F}_q.$ 
			\item[b)] If $p\notdivides  w$, then linear code spanned by the rows of the matrix $A=[O, \sqrt{-wd}\textbf{1}]$ is a self-orthogonal $[m+1,\rm{rank} O]$ code over the field $\mathbb{F},$ where $\mathbb{F}=\mathbb{F}_q$ if $-wd$ is square root in $\mathbb{F}_q,$ and $\mathbb{F}=\mathbb{F}_{q^2}$ otherwise.
	\end{itemize}

	\item If $a\neq d,$ we differ three cases.
	\begin{itemize}
			\item[a)] If $p\mid w$, then linear code spanned by the rows of the matrix $A=[\sqrt{d-a} I_m, O]$ is a self-orthogonal $[m+n,m]$ code over the field $\mathbb{F},$ where $\mathbb{F}=\mathbb{F}_q$ if $d-a$ is square root in $\mathbb{F}_q,$ and $\mathbb{F}=\mathbb{F}_{q^2}$ otherwise. If $m=n,$ the obtained code is self-dual.
			\item[b)] If $p\mid w-1$, then linear code spanned by the rows of the matrix $A=[\sqrt{wd-a} I_m, O, \sqrt{-wd}\textbf{1}]$ is a self-orthogonal $[m+n+1,m]$ code over the field $\mathbb{F},$ where $\mathbb{F}=\mathbb{F}_q$ if $wd-a$ and $-wd$ are square roots in $\mathbb{F}_q,$ and $\mathbb{F}=\mathbb{F}_{q^2}$ otherwise.
			\item[c)] If $p\notdivides  w$ and $p\notdivides  w-1$, then binary linear code spanned by the rows of the matrix $A=[\sqrt{d-a} I_m, O, \sqrt{-wd}\textbf{1}]$ is a self-orthogonal $[m+n+1,m]$ code over the field $\mathbb{F},$ where $\mathbb{F}=\mathbb{F}_q$ if $d-a$ and $-wd$ are square roots in $\mathbb{F}_q,$ and $\mathbb{F}=\mathbb{F}_{q^2}$ otherwise.
	\end{itemize}
\end{itemize}
\end{thm}
\textbf{Proof:}
\begin{itemize}
		\item If $a=d,$ it follows form Remark \ref{remarkmodp} that $O[s]\cdot O[t] \jmod{dw}{p},$ for $s, t\in \{1,\ldots,m\}.$ It is now easy to see that $O$ generates self-orthogonal code if $p\mid w,$ and that $A=[O,\sqrt{-dw}\textbf{1}]$ generates self-orthogonal code otherwise. 
		\item If $a\neq d,$ it follows from Remark \ref{remarkmodp} that $O[s]\cdot O[t] \jmod{dw}{p},$ for $s\neq t$ and $ O[s]\cdot O[s] \jmod{a-d+wd}{p}.$
				\begin{itemize}
							\item[a)] If $p\mid w,$ then $O[s]\cdot O[t]\jmod{0}{p},$ and $O[s]\cdot O[t]\jmod{a-d}{p}.$ The matrix $[\sqrt{d-a}I_m, O]$ generates self-orthogonal 
								$[m+n,m]$ code. If $m=n,$ dimension of the code is half of its length, so the code is self-dual.
							\item[b)] If $p\mid w-1,$ then $O[s]\cdot O[t]\jmod{wd}{p},$ and $O[s]\cdot O[t]\jmod{a}{p}.$ The matrix $[\sqrt{wd-a}I_m, O, \sqrt{-wd}\textbf{1}0]$ generates self-orthogonal $[m+n+1,m]$ code.
							\item[c)] If $p\notdivides  w$ and $p\notdivides  w-1$, then $O[s]\cdot O[t]\jmod{wd}{p},$ and $O[s]\cdot O[t]\jmod{a-d+wd}{p}.$ The matrix $[\sqrt{d-a}I_m, O, \sqrt{-wd}\textbf{1}0]$ generates self-orthogonal $[m+n+1,m]$ code.
				\end{itemize}
\end{itemize}
\begin{flushright}
$\blacksquare$
\end{flushright}

%

\begin{thm}\label{ft1}
Let $\mathcal{D}$ be a weakly self-orthogonal $1-(v,k,r)$ design such that $k$ is odd and block intersection numbers are odd. Let $G$ be an automorphism group of $\mathcal{D}$ acting with $f_1$ fixed points and $n$ point orbits of length $2,$ and with $f_2$ fixed blocks and $m$ block orbits of length $2.$
Then
\begin{enumerate}
	\item the binary linear code spanned by the matrix $[OM1,\textbf{1}]$ is a self-orthogonal code of length $f_1+1$ and dimension $\rm{rank}(OM1)$,
	\item the binary linear code spanned by the matrix $ OM2$ is a self-orthogonal code of length $n$ and dimension $\rm{rank}(OM2).$
\end{enumerate}
\end{thm}
\textbf{Proof:}
 \begin{enumerate}
	\item Since $k$ is odd, each block contains odd number of fixed points and since block intersection numbers are odd, intersection of every two blocks contains odd number of fixed points. Therefore, it follows that the matrix $[OM1,\textbf{1}]$ spans a binary self-orthogonal code.

	\item Let $\mathcal{B}_t$ and $\mathcal{B}_s$ be block orbits of size $2$ under the action of the group $G$ on the design.
	For $s\neq t,$ since $\mathcal{B}_t$ is block orbit of size 2, it follows from Remark \ref{remark_sum} that 
	$$\sum_{x'\in\mathcal{B}_t} |x\cap x'| \jmod{2}{2} \jmod{0}{2},$$
	and for $s=t,$ it follows that
	$$\sum_{x'\in\mathcal{B}_s} |x\cap x'|=|x\cap x|+ \sum_{x'\in\mathcal{B}_s, \\ x\neq x'}|x\cap x'| \jmod{1+1}{2}.$$
	We conclude that 
	$\sum_{j=f_1+1}^{f_1+n}a_{s,j}a_{t,j} \jmod{0}{2},$ for $s\neq t$ and $\sum_{j=f_1+1}^{f_1+n}a_{s,j}a_{s,j} \jmod{0}{2}$ and that the binary code spanned by the matrix $OM2$ is self-orthogonal.
\end{enumerate}
\begin{flushright}
$\blacksquare$
\end{flushright}
We generalize previous theorem to obtain self-orthogonal codes over $\mathbb{F}_q,$ where $q$ is a prime power.

\begin{thm}\label{ftsl4_q}
Let $q=p^l$ be prime power and $\mathbb{F}_q$ a finite field of order $q.$ Let $\mathcal{D}$ be $1-(v,k,r)$ design such that $k \jmod{a}{p}$ and $|B_i\cap B_j| \jmod{d}{p},$ for all $i,j\in\{1,\ldots,b\}, \ i\neq j,$ where $B_i$ and $B_j$ are two blocks of the design $\mathcal{D},$ and let $G$ be an automorphism group of the design which acts on $\mathcal{D}$ with $f_1$ fixed points and $n$ point orbits of length $p^\alpha, \ 1\leq\alpha\leq l$, and with $f_2$ fixed blocks and $m$ block orbits of length $p^\alpha$. 

Then 
\begin{itemize}
	\item If $a=q,$ we differ two cases. 
			\begin{enumerate}
					\item the linear code spanned by the matrix $[OM1,\sqrt{-a}\textbf{1}]$ is a self-orthogonal $[f_1+1,\rm{rank}(OM1)]$ code over the field $\mathbb{F},$
					where $\mathbb{F}=\mathbb{F}_q$ if $-a$ is square root in $\mathbb{F}_q,$ and $\mathbb{F}=\mathbb{F}_{q^2}$ otherwise.
					\item the linear code spanned by the matrix $OM2$ is a self-orthogonal $[n,\rm{rank}(OM2)]$ code of length over the field $\mathbb{F}_q.$
			\end{enumerate}
	
	\item If $a\neq d,$ we differ two cases.
			\begin{enumerate}
					\item the linear code spanned by the matrix $[\sqrt{d-a}\cdot I_{f_2}, OM1,\sqrt{-d}\textbf{1}]$ is a self-orthogonal $[f_2+f_1+1,f_2]$ code over 
					the field $\mathbb{F},$ where $\mathbb{F}=\mathbb{F}_q$ if $d-a$ and $-d$ are square roots in $\mathbb{F}_q,$ and $\mathbb{F}=\mathbb{F}_{q^2}$ otherwise.
					\item the linear code spanned by the matrix $[\sqrt{d-a}\cdot I_m, OM2]$ is a self-orthogonal $[m+n,m]$ code over the field $\mathbb{F},$ where 
					$\mathbb{F}=\mathbb{F}_q$ if $d-a$ is square root in $\mathbb{F}_q,$ and $\mathbb{F}=\mathbb{F}_{q^2}$ otherwise. If $m=n,$ the obtained code is self-dual.
			\end{enumerate}
	\end{itemize}
\end{thm}
\textbf{Proof:}
 \begin{enumerate}
	\item Since $k\jmod{a}{p}$, each block contains $q\cdot \beta+a$ fixed points and since $|B_i\cap B_j| \jmod{d}{p},$ for all $i,j\in\{1,\ldots,b\}, \ i\neq j,$ intersection of every two blocks contains $q\cdot \gamma+d$ fixed points. 
	Therefore, if $a=d$, it follows that the matrix $[OM1,\sqrt{-a}\textbf{1}]$ spans a self-orthogonal $[f_1+1,\rm{rank}(OM1)]$ code over the field $\mathbb{F},$ and if $a\neq d$, it follows that the matrix $[\sqrt{d-a}\cdot I_{f_2}, OM1,\sqrt{-d}\textbf{1}]$ spans a self-orthogonal $[f_2+f_1+1,f_2]$ code over the field $\mathbb{F}.$

	\item For $s\neq t,$ since $\mathcal{B}_t$ is block orbit of size $p^\alpha$, it follows from Remark \ref{remark_sum} that 
	$$\sum_{x'\in\mathcal{B}_t} |x\cap x'| \jmod{p^\alpha d}{p} \jmod{0}{p},$$
	and for $s=t$ it follows that
	$$\sum_{x'\in\mathcal{B}_s} |x\cap x'|=|x\cap x|+ \sum_{x'\in\mathcal{B}_s, \\ x\neq x'}|x\cap x'| \jmod{a+(p^\alpha-1)d}{p}.$$
	It follows that 
	$OM2[s]\cdot OM2[t] \jmod{0}{p},$ for $s\neq t$ and $OM2[s]\cdot OM2[s] \jmod{a-d}{p}.$ If $a=d$, then the code spanned by the matrix $OM2$ is self-orthogonal $[n,\rm{rank}(OM2)]$ code over $\mathbb{F}_q$, and if $a\neq d,$ it follows that the code spanned by the matrix $[\sqrt{d-a}\cdot I_m, OM2]$ is self-orthogonal $[nm+n,n]$ code over the field $\mathbb{F}.$ If $m=n,$ dimension of the code is half of its length, so the code is self-dual.
\end{enumerate}
\begin{flushright}
$\blacksquare$
\end{flushright}

\section{Results}\label{results}

\begin{thm}
	Let $G$ be a finite permutation group acting transitively on the set $\Omega$ of size $n.$ Let $\alpha\in\Omega$ and $\Delta=\bigcup_{i=1}^s \delta_iG_\alpha,$ where $\delta_i,\ldots,\delta_s\in\Omega$ are representatives of distinct $G_\alpha$-orbits. If $\Delta\neq\Omega$ and $$\mathcal{B}=\{\Delta g\mid g\in G\},$$ then $\mathcal{D}=(\Omega,\mathcal{B})$ is $1-(n,|\Delta|,\frac{|G_\alpha|}{|G_{\Delta}|}\sum_{i=1}^n |\alpha G_{\delta_i}|)$ design with $\frac{m\cdot|G_\alpha|}{|G_{\Delta}|}$ blocks.
\end{thm}
\textbf{Proof:} {} In \cite{1}.
\begin{flushright}
	$\blacksquare$
\end{flushright}
\noindent
Using this construction, we constructed 178 pairwise non-isomorphic weakly self-orthogonal 1-designs on 22, 55, 66, 110, 132, 144 and 165 points invariant under an action of $M_{11}.$ Precisely, 4 designs on 22 points, 2 designs on 55 points, 6 designs on 66 points, 41 designs on 110 points, 76 designs on 132 points, 26 designs on 144 points and 20 designs on 165 points.
\\
Using Theorem \ref{tmTonchev}, from constructed designs, we obtained binary self-orthogonal codes. Optimal codes will be denoted with $*$, best known codes will be denoted with $(**)$, and near optimal codes will be denoted with $+$. 
All tables will be ordered by four cases of weakly self-orthogonal designs.
\begin{enumerate}
	\item[Case 1.] Codes obtained from self-orthogonal designs.
	\item[Case 2.] Codes obtained from weakly self-orthogonal designs such that block sizes are even and block intersection numbers are odd.
	\item[Case 3.] Codes obtained from weakly self-orthogonal designs such that block sizes are odd and block intersection numbers are even.
	\item[Case 4.] Codes obtained from weakly self-orthogonal designs such that block sizes are odd and block intersection numbers are odd.
\end{enumerate}

\begin{table}[H]
	\caption{Non-trivial binary pairwise non-equivalent self-orthogonal codes obtained from non-trivial self-orthogonal designs using Theorem \ref{tmTonchev} (Case 1)}
	\scriptsize 
	\begin{tabular}{lll}
		\hline\noalign{\smallskip}
		\textbf{$C$}& \textbf{Aut$C$} or \textbf{$\#$Aut$C$} & \textbf{Design} \\ 
		\noalign{\smallskip}\hline\noalign{\smallskip}
		$[22,10,4]$ 	 & $Z_2\times(E_{2^10}:S_{11})$ & $1-(22,10,10)$ \\
		$[22,11,2]$    & $Z_2\times(E_{2^10}:S_{11})$ & $1-(22,2,1)$ \\\hline
		$[66,10,20]$   & $S_{12}$ 				 & $1-(66,20,20)$ \\
		$[66,11,20]$   & $S_{12}$ 			 	& $1-(66,46,46)$ \\\hline
		$[110,44,8]$   &  $2^{55}\cdot 11!$  		& $1-(110,72,36)$ \\
		$[110,10,20]$  &  $2^{55}\cdot 11!$ 	  & $1-(110,36,18)$ \\
		$[110,54,4]$   &  $2^{55}\cdot 55!$  	  & $1-(110,108,54)$ \\
		$[110,11,20]$  &  $2^{55}\cdot 11!$   	& $1-(110,74,37)$ \\
		$[110,45,6]$   &   $2^{55}\cdot 11!$ 	  & $1-(110,38,19)$ \\
		$[110,11,10]$  &   $11!\cdot (10!)^{11}$ & $1-(110,10,1)$ \\
		$[110,55,6]$   &   $S_{11}$							& $1-(110,18,9)$ \\
		$[110,11,20]$  &   $2^{55}\cdot 11!$ 		& $1-(110,74,37)$ \\
		$[110,10,20]$  &  $11!\cdot (10!)^{11}$ & $1-(110,20,10)$ \\
		$[110,54,8]$   & $S_{11}$								& $1-(110,28,28)$ \\\hline
		$[132,54,16]$  & $M_{11}$  							& $1-(132,20,20)$ \\
		$[132,22,6]$   &  $22!\cdot(6!)^{22}$   & $1-(132,6,1)$ \\
		$[132,21,12]$  &  $22!\cdot(6!)^{22}$   & $1-(132,60,60)$ \\
		$[132,65,12]$  & $M_{11}$	   				    & $1-(132,40,40)$ \\
		$[132,66,6]$   & $M_{11}$						    & $1-(132,26,26)$ \\
		$[132,11,12]$  &  $11!\cdot(12!)^{11}$  & $1-(132,12,1)$ \\
		$[132,11,24]$  &  $61440\cdot 12!\cdot(6!)^{21}$   & $1-(132,66,6)$ \\
		$[132,55,16]$  & $M_{11}$							    & $1-(132,46,46)$ \\
		$[132,55,16]$  & $M_{11}$						    & $1-(132,46,46)$ \\
		$[132,45, 20]$ & $M_{11}$ 							& $1-(132,32,32)$ \\
		$[132,44,24]$  & $M_{11}$						    & $1-(132,100,100)$ \\
		$[132,10,24]$  & $11!\cdot(12!)^{11}$  & $1-(132,120,10)$ \\
		$[132,55,12]$  & $M_{11}$							    & $1-(132,112,112)$ \\
		$[132,54,8]$   &    $2^{66}\cdot 12!$			& $1-(132,60,30)$ \\
		$[132,22,22]$  &    $2^{66}\cdot 7920$		& $1-(132,30,15)$ \\
		$[132,66,2]$   &    $2^{66}\cdot 66!$			& $1-(132,2,1)$ \\
		$[132,65,4]$   &    $2^{66}\cdot 66!$			& $1-(132,20,20)$ \\
		$[132,45,16]$  &    $2^{66}\cdot 7920$		& $1-(132,32,16)$ \\
		$[132,21,32]$  &    $2^{66}\cdot 7920$		& $1-(132,80,80)$ \\
		$[132,55,6]$   &    $2^{66}\cdot 12!$			& $1-(132,50,50)$ \\
		$[132,11,40]$  &    $2^{66}\cdot 12!$			& $1-(132,92,46)$ \\
		$[132,11,22]$  &    $2^{66}\cdot 12!$			& $1-(132,22,2)$ \\
		$[132,11,40]$  &    $2^{66}\cdot 12!$			& $1-(132,110,10)$ \\
		$[132,10,40]$  &    $2^{66}\cdot 12!$			& $1-(132,40,20)$ \\
		$[132,55,8]$   &    $2^{66}\cdot 12!$			& $1-(132,82,82)$ \\
		$[132,44,16]$  &    $2^{66}\cdot 7920$		& $1-(132,100,50)$ \\
		$[132,55,8]$   &    $2^{66}\cdot 12!$		  & $1-(132,72,36)$ \\\hline
		$[144,12,12]$  &    $945^{13}\cdot 7920^{13}\cdot 2^{78}$& $1-(144,12,1)$ \\
		$[144,56,20]$  &$M_{12}:Z_2$			        & $1-(144,56,56)$ \\
		$[144,46,22]$  &$M_{12}:Z_2$		        & $1-(144,66,30)$ \\\hline
		$[165,54,20]$  &$M_{11}$ 		          & $1-(165,48,48)$ \\
		$[165,44,32]$  &$M_{11}$		          & $1-(165,72,72)$ \\
		$[165,10,72]$  &$S_{11}$  			            & $1-(165,80,80)$ \\
		\noalign{\smallskip}\hline
	\end{tabular}
\end{table}

\begin{table}[H]
	\caption{Non-trivial binary self-orthogonal codes obtained from non-trivial weakly self-orthogonal designs using Theorem \ref{tmTonchev} (Case 2)}
	\footnotesize
	\begin{tabular}{lll}
		\hline\noalign{\smallskip}
		\textbf{$C$}& \textbf{Aut$C$} or \textbf{$\#$Aut$C$} & \textbf{Design} \\ 
		\noalign{\smallskip}\hline\noalign{\smallskip}
		$[331,165,12]$  &$$  			            & $1-(165,116,116)$ \\
		$[331,165,8]$   &$$  			            & $1-(165,84,84)$ \\
		$[331,165,12]$  &$ $ 			        & $1-(165,36,36)$ \\
		\noalign{\smallskip}\hline
	\end{tabular}
\end{table}

\begin{table}[H]
	\caption{Non-trivial binary self-orthogonal codes obtained from non-trivial weakly self-orthogonal designs using Theorem \ref{tmTonchev} (Case 3)}
	\scriptsize
	\begin{tabular}{lll}
		\hline\noalign{\smallskip}
		\textbf{$C$}& \textbf{Aut$C$} or \textbf{$\#$Aut$C$} & \textbf{Design} \\ 
		\noalign{\smallskip}\hline\noalign{\smallskip}
		$[132,66,6]$ & 			 $						$	& $1-(66,21,21)$ \\
		$[132,66,8]$ & 			 $						$	& $1-(66,45,45)$ \\
		$[220,110,4]$ & 		 $						$	  & $1-(110,73,73))$ \\
		$[220,110,4]$ & 		 $					$	& $1-(110,37,37))$ \\
		$[220,110,4]$ & 		 $						$  	& $1-(110,73,73)$  \\
		$[220,110,4]$ & 		 $						$  	& $1-(110,37,37)$ \\
		$[220,110,4]$ & 		 $						$  	& $1-(110,9,9)$ \\
		$[220,110,4]$ & 		 $						$  	& $1-(110,73,73) $ \\
		$[220,110,8]$ & 		 $						$  	& $1-(110,81,81)$ \\
		$[220,110,8]$ & 		 $						$  	& $1-(110,29,29)$ \\
		$[220,110,4]$ & 		 $						$  	& $1-(110,37,37)$ \\
		$[220,110,4]$ & 		 $						$  	& $1-(110,101,101)$ \\\hline
		$[264,132,4]$ &			 $							$		& $1-(132,5,5)$ \\
		$[144,12,12]$ &			 $							$		& $1-(132,11,1)$ \\
		$[264,132,12]$ &		 $							$	& $1-(132,21,21)$ \\
		$[264,132,4]$ &			 $							$		& $1-(132,7,7)$ \\
		$[264,132,10]$ &		 $							$	& $1-(132,25,25)$ \\
		$[264,132,4]$ &			 $							$	& $1-(132,11,11)$ \\
		$[264,132,12]$ &		 $  					$	& $1-(132,27,27)$ \\
		$[144,12,56]$ &			 $					  	$	& $1-(132,55,5)$ \\
		$[264,132,12]$ &		 $							$	& $1-(132,31,31)$ \\
		$[264,132,12]$ &			 $							$	& $1-(132,101,101)$ \\
		$[144,12,58]$ &			 $						$	& $1-(132,77,7)$ \\
		$[264,132,12]$ &			 $							$	& $1-(132,105,105)$ \\
		$[264,132,4]$ &			 $							$	& $1-(132,121,121)$ \\
		$[264,132,12]$ &			 $							$	& $1-(132,107,107)$ \\
		$[264,132,4]$ &			 $							$	& $1-(132,125,125)$ \\
		$[264,132,12]$ &			 $							$	& $1-(132,111,111)$ \\
		$[144,12,22]$ &			 $							$	& $1-(132,121,11)$ \\
		$[264,132,4]$ &			 $							$	& $1-(132,127,127)$ \\
		$[264,132,4]$ &			 $							$	& $1-(132,61,61)$ \\
		$[144,12,12]$ &			 $							$	& $1-(132,11,1)$ \\
		$[264,132,4]$ &			 $							$	& $1-(132,31,31)$ \\
		$[264,132,4]$ &			 $							$	& $1-(132,91,91)$ \\
		$[264,132,4]$ &			 $							$	& $1-(132,41,41)$ \\
		$[264,132,4]$ &			$							$	& $1-(132,101,101)$ \\
		$[144,12,22]$ &			 $							  $	& $1-(132,121,11) $ \\
		$[264,132,4]$ &			 $						$	& $1-(132,71,71) $ \\\hline
		$[288,144,4]$ & 	& $1-(144,11,11)$ \\
		$[288,144,12]$ &	& $1-(144,11,11)$ \\
		$[288,144]$	&& $1-(144,55,55)$ \\
		$[288,144]$	&& $	1-(144,23,23)$ \\
		$[288,144]	$	&& $1-(144,23,23)$ \\
		$[288,144]$	&& $	1-(144,67,67)$ \\
		$[288,144]$	&& $	1-(144,67,67)$ \\
		$[288,144]$	&& $	1-(144,67,67)$ \\
		$[288,144]$	&& $	1-(144,67,67)$ \\
		$[288,144]$	&& $	1-(144,77,77)$ \\
		$[288,144]$	&& $	1-(144,77,77)$ \\
		$[288,144]$	&& $	1-(144,77,77)$ \\
		$[288,144]$	&& $	1-(144,77,77)$ \\
		$[288,144]$	&& $	1-(144,121,121)$ \\
		$[288,144]$	&& $	1-(144,121,121)$ \\
		$[288,144]$	&& $	1-(144,89,89)$ \\
		$[288,144]$	&& $	1-(144,133,133)$ \\
		$[288,144]$	&& $	1-(144,133,133)$ \\ \hline
		$[330,165,12]$  &$$ 		    & $1-(165,129,129)$ \\
		$[330,165,8]$  &$$  	& $1-(165,81,81)$ \\
		$[330,165,12]$  &$ $ & $1-(165,49,49)$ \\
		\noalign{\smallskip}\hline
	\end{tabular}
\end{table}

\begin{table}[H]
	\caption{Non-trivial binary pairwise non-equivalent self-orthogonal codes obtained from non-trivial weakly self-orthogonal designs using Theorem \ref{tmTonchev} (Case 4)}
	\footnotesize
	\begin{tabular}{lll}
		\hline\noalign{\smallskip}
		\textbf{$C$}& \textbf{Aut$C$} or \textbf{$\#$Aut$C$} & \textbf{Design} \\ 
		\noalign{\smallskip}\hline\noalign{\smallskip}
		$[166,55,20]$  &$M_{11}$ 			    & $1-(165,109,109)$ \\
		$[166,45,26]$  &$M_{11}$  			    & $1-(165,61,61)$ \\
		$[166,11,46]$  &$S_{11}$ 		        & $1-(165,85,85)$ \\
		\noalign{\smallskip}\hline
	\end{tabular}
\end{table}

\noindent
Using Theorems \ref{tmsl1},\ref{tmsl2},\ref{tmsl3} and \ref{tmsl4}, we constructed binary self-orthogonal codes from the orbit matrices of the non-trivial weakly self-orthogonal 1-designs on less than 165 points (including). In order to construct the orbit matrices we determined orbits of all cyclic subgroups of prime order of the group $M_{11}$ acting with orbits of the same length. 

\begin{table}[H]
	\caption{Non-trivial binary pairwise non-equivalent self-orthogonal codes obtained from the orbit matrices of the non-trivial self-orthogonal 1-designs on 110 points using Theorem \ref{tmsl1} (Case 1)}
	\label{t1}
	\footnotesize
	\begin{tabular}{lll}
		\hline\noalign{\smallskip}
		\textbf{$C$}& \textbf{Aut$C$} or \textbf{$\#$Aut$C$} & \textbf{Design} \\ 
		\noalign{\smallskip}\hline\noalign{\smallskip}
		$[22,8,4]$ & $Z_2\times((((((Z_2\times D_8):Z_2):Z_3):Z_2):Z_2)\times((E_{2^4}:A_5):Z_2)\times S_4)$	& \multirow{2}{*}{$1-(110,72,36)$}\\ 
		$[10,4,4]$* &   $Z_2\times((E_{2^4}:A_5):Z_2)$								 										&  	\\
		$[22,2,4]$ & $S_{10}\times S_8\times D_8$ 														& $1-(110,36,18)$ \\
		$[10,5,2]$ &   $Z_2\times((E_{2^4}:A_5):Z_2)$							 										& $1-(110,2,1)$ \\
		$[22,3,4]$ & $S_{10}\times S_8\times D_8$ 										& $1-(110,74,37)$\\
		$[22,9,2]$ & $Z_2\times((((((Z_2\times D_8):Z_2):Z_3):Z_2):Z_2)\times((E_{2^4}:A_5):Z_2)\times S_4)$ &  $1-(110,38,19)$\\
		$[22,3,2]$ &  $2^{18}\cdot 3^8\cdot 5^4\cdot 7^2$		& $1-(110,10,1)$  \\
		$[22,11,2]$ & $Z_2\times(((E_{2^5}:A_6):Z_2)\times((E_{2^4}:A_5):Z_2)$		& $1-(110,18,9)$  \\ 
		$[22,2,12]$ & $2^{18}\cdot 3^8\cdot 5^4\cdot 7^2$		& $1-(110,20,10)$  \\
		\noalign{\smallskip}\hline
	\end{tabular}
\end{table}

\begin{table}[H]
	\caption{Non-trivial binary pairwise non-equivalent self-orthogonal codes obtained from the orbit matrices of the non-trivial self-orthogonal 1-designs on 132 points using Theorem \ref{tmsl1} (Case 1)}
	\label{t2}
	\footnotesize
	\begin{tabular}{lll}
		\hline\noalign{\smallskip}
		\textbf{$C$}& \textbf{Aut$C$} or \textbf{$\#$Aut$C$} & \textbf{Design} \\ 
		\noalign{\smallskip}\hline\noalign{\smallskip}
		$[12,4,4]$ & $E_{2^2}\times((E_{2^4}:A_5):Z_2)$ & $1-(132,20,20)$  \\
		$[12,2,6]$ & $(S_6\times S_6):Z_2$ 		& $1-(132,6,1)$  \\ 
		$[12,5,4]$* & $((E_{2^5}:A_6):Z_2):Z_2$ 	& $1-(132,40,40)$  \\
		$[12,5,4]$* & $(E_{2^4}:A_5):Z_2$ 		& $1-(132,46,46)$  \\
		$[12,2,2]$ & $Z_2\times S_{10}$  &$1-(132,30,15)$  \\
		$[12,5,2]$ & $E_{2^2}\times((E_{2^4}:A_5):Z_2)$ & $1-(132,50,50)$  \\
		$[12,11,2]$* &  $2^9\cdot 3^4\cdot 5^2\cdot 7$		& $1-(132,22,2)$  \\ 
		$[12,1,10]$ &  $2^9\cdot 3^4\cdot 5^2\cdot 7$		& $1-(132,11,40)$ \\
		$[12,5,2]$ & $E_{2^2}\times((E_{2^4}:A_5):Z_2)$ 		& $1-(132,82,82)$  \\
		\noalign{\smallskip}\hline
	\end{tabular}
\end{table}

\begin{table}[H]
	\caption{Non-trivial binary pairwise non-equivalent self-orthogonal codes obtained from the orbit matrices of the non-trivial self-orthogonal 1-designs on 144 points using Theorem \ref{tmsl1} (Case 1)}
	\label{t3}
	\footnotesize
	\begin{tabular}{lll}
		\hline\noalign{\smallskip}
		\textbf{$C$}& \textbf{Aut$C$} or \textbf{$\#$Aut$C$} & \textbf{Design} \\ 
		\noalign{\smallskip}\hline\noalign{\smallskip}
		$[72,4,12]$ &  $7920^6\cdot 2^{41}\cdot 3^{19}\cdot 5^6\cdot 7^7\cdot 13\cdot 17\cdot 19\cdot 23$ & \multirow{2}{*}{$1-(144,12,1)$} \\ 
		$[48,6,4]$ &  $2^{41}\cdot 3^{20}\cdot 5^6\cdot 7^3\cdot 11^3$ & \\
		$[72,4,12]$ & $((E_{2^4}:Z_3):Z_2):Z_2$ &\multirow{2}{*}{ $1-(144,56,56)$} \\ 
		$[48,20,8]$ & $Z_2\times (S_4\times S_4):Z_2)$ &  \\
		$[72,22,12]$ & $((E_{2^4}:Z_3):Z_2):Z_2$  &\multirow{2}{*}{ $1-(144,66,30)$} \\ 
		$[48,16,8]$ & $(S_3\times S_3):Z_2$ & \\
		\noalign{\smallskip}\hline
	\end{tabular}
\end{table}

\begin{table}[H]
	\caption{Non-trivial binary pairwise non-equivalent self-orthogonal codes obtained from the orbit matrices of the non-trivial self-orthogonal 1-designs on 165 points using Theorem \ref{tmsl1} (Case 1)}
	\label{t4}
	\footnotesize
	\begin{tabular}{lll}
		\hline\noalign{\smallskip}
		\textbf{$C$}& \textbf{Aut$C$} or \textbf{$\#$Aut$C$} & \textbf{Design} \\ 
		\noalign{\smallskip}\hline\noalign{\smallskip}
		$[33,10,4]$& $2048$ & \multirow{2}{*}{$1-(165,48,48)$} \\ 
		$[15,4,8]$*  & $A_8$  & \\
		$[33,8,8]$  &$E_{2^2}\times ((E_{2^5}\times D_8):Z_2)$ & $1-(165,72,72)$\\
		$[33,2,16]$  & $2^{27}\cdot 3^{12}\cdot 5^5\cdot 7^2\cdot 11^2$ & $1-(165,80,80)$ \\
		\noalign{\smallskip}\hline
	\end{tabular}
\end{table}

\begin{table}[H]
	\caption{Non-trivial binary pairwise non-equivalent self-orthogonal codes obtained from the orbit matrices of the non-trivial weakly self-orthogonal 1-designs on 165 points using Theorem \ref{tmsl2} (Case 2)}
	\label{tab1}
	\footnotesize
	\begin{tabular}{lll}
		\hline\noalign{\smallskip}
		\textbf{$C$}& \textbf{Aut$C$} or \textbf{$\#$Aut$C$} & \textbf{Design} \\ 
		\noalign{\smallskip}\hline\noalign{\smallskip}
		$[67,33,4]$ &	$2^{14} $			& \multirow{2}{*}{$1-(165,116,116)$}\\  
		$[31,15,8]$* & $ \rm{PSL}(5,2)$  					&  	\\
		$[67,33,4]$ & $  $ 			& \multirow{2}{*}{$1-(165,84,84)$}\\  
		$[31,15,4]$ &  $2^{25}\cdot 3^6\cdot 5^3\cdot 7^2\cdot 11\cdot 13$	 					&  	\\
		$[67,33,4]$ & $E_{2^3}\times ((Z_2\times(Z_4\times Z_2):Z_2):Z_2)):Z_2)$ & $1-(165,36,36)$ \\
		\noalign{\smallskip}\hline
	\end{tabular}
\end{table}

\begin{table}[H]
	\caption{Non-trivial binary pairwise non-equivalent self-orthogonal codes obtained from the orbit matrices of the non-trivial weakly self-orthogonal 1-designs on 22, 55 and 66 points using Theorem \ref{tmsl3} (Case 3)}
	\label{sl3t0}
	\footnotesize
	\begin{tabular}{lll}
		\hline\noalign{\smallskip}
		\textbf{$C$}& \textbf{Aut$C$} or \textbf{$\#$Aut$C$} & \textbf{Design} \\ 
		\noalign{\smallskip}\hline\noalign{\smallskip}
		$[12,6,2]$ & $((E_{2^5}:A_6):Z_2):Z_2$ & $1-(66,21,21)$ \\
		$[12,6,4]$* & $(E_{2^5}:A_6):Z_2)$   & $1-(66,21,21)$ \\
		\noalign{\smallskip}\hline
	\end{tabular}
\end{table}

\begin{table}[H]
	\caption{Non-trivial binary pairwise non-equivalent self-orthogonal codes obtained from the orbit matrices of the non-trivial weakly self-orthogonal 1-designs on 110 points using Theorem \ref{tmsl3} (Case 3)}
	\label{sl3t1}
	\footnotesize
	\begin{tabular}{lll}
		\hline\noalign{\smallskip}
		\textbf{$C$}& \textbf{Aut$C$} or \textbf{$\#$Aut$C$} & \textbf{Design} \\ 
		\noalign{\smallskip}\hline\noalign{\smallskip}
		$[44,22,4]$ & 	$2^{19}\cdot 10! \cdot 8! \cdot 4!$			& $1-(110,73,73)$\\
		$[44,22,2]$ & 	$2^{20}\cdot 10! \cdot 8! \cdot 4! $& $1-(110,37,37)$ \\
		$[44,22,2]$ & 	$2^{20} \cdot 10!\cdot 15!\cdot 8!\cdot 4! $& $1-(110,9,9)$ \\
		$[44,22,4]$ & $2^{30}\cdot 3^4 \cdot 5^2$& $1-(110,81,81)$ \\
		$[44,22,4]$ & 	$2^{37}\cdot 3^8 \cdot 5^4\cdot 7^2$& $1-(110,101,101)$ \\
		\noalign{\smallskip}\hline
	\end{tabular}
\end{table}

\begin{table}[H]
	\caption{Non-trivial binary pairwise non-equivalent self-orthogonal codes obtained from the orbit matrices of the non-trivial weakly self-orthogonal 1-designs on 132 points using Theorem \ref{tmsl3} (Case 3)}
	\label{sl3t2}
	\footnotesize
	\begin{tabular}{lll}
		\hline\noalign{\smallskip}
		\textbf{$C$}& \textbf{Aut$C$} or \textbf{$\#$Aut$C$} & \textbf{Design} \\ 
		\noalign{\smallskip}\hline\noalign{\smallskip}
		$[24,12,4]$ & 	$E_{2^{10}}:(A_6\times A_6):D_8$ & $1-(132,5,5)$ \\
		$[14,2,2]$ & 	$Z_2\times S_{12} $& $1-(132,11,1)$ \\
		$[24,12,2]$ & 	$D_8 x (E_{2^{10}} : S_5)$& $1-(132,21,21)$ \\
		$[24,12,4]$ & 	$E_{2^{10}}:(A_6\times A_6):D_8$ & $1-(132,7,7)$ \\
		$[24,12,6]$ & $((E_{2^6}:(Z_2.A_6)):Z_2 $ & $1-(132,25,25)$ \\
		$[24,12,4]$ & 	$E_{2^11}\cdot S_{12} $& $1-(132,11,11)$ \\
		$[24,12,8]*$ & $M_{24} $ & $1-(132,27,27)$ \\
		$[14,2,6]$ & 	$S_6\times S_8$ & $1-(132,55,5)$ \\
		$[24,12,4]$ & 	$E_{2^{10}} : ((Z_3 . A_6) : Z_2)$& $1-(132,31,31)$ \\
		$[24,12,2]$ & 	$2^{12}\cdot 10! $& $1-(132,61,61)$ \\
		\noalign{\smallskip}\hline
	\end{tabular}
\end{table}

\begin{table}[H]
	\caption{Non-trivial binary self-orthogonal codes obtained from the orbit matrices of the non-trivial weakly self-orthogonal 1-designs on 144 points using Theorem \ref{tmsl3} (Case 3)}
	\label{sl3t3}
	\footnotesize
	\begin{tabular}{lll}
		\hline\noalign{\smallskip}
		\textbf{$C$}& \textbf{Aut$C$} or \textbf{$\#$Aut$C$} & \textbf{Design} \\ 
		\noalign{\smallskip}\hline\noalign{\smallskip}
		$[144,72,2]$ & 	$ $			& \multirow{2}{*}{$1-(144,11,11)$}\\  
		$[96,48,4]$ &  $2^{83}\cdot 3^20\cdot 7^6\cdot 5^6 \cdot 11^3$	 					&  	\\
		$[144,72,10]$  &	$ $			& \multirow{2}{*}{$1-(144,11,11)$}\\  
		$[96,48,12]$ & $(S_2\times S_3):Z_2 $ 	 					&  	\\
		$[144,72,12]$ & $ $			& \multirow{2}{*}{$1-(144,55,55)$}\\  
		$[96,48,8]$ & $Z_2\times Z_3\times S_3 $ 	 					&  	\\
		$[144,72,12]$ & $ $ 	& \multirow{2}{*}{$1-(144,23,23)$}\\  
		$[96,48,16]**$ & $D_{12} $ 	 					&  	\\
		$[144,72,8]$ & $ $ 		& \multirow{2}{*}{$1-(144,67,67)$}\\  
		$[96,48,8]$ & $Z_2\times ((S_3\times S_3):Z_2) $ 	 					&  	\\
		$[144,72,6]$ & $ $ 			& \multirow{2}{*}{$1-(144,67,67)$}\\  
		$[96,48,4]$ & $$  					&  	\\
		$[96,48,8]$ & $(E_{2^3}\times S_3):Z_2$& $1-(144,67,67)$ \\
		$[144,72,6]$&  &\multirow{2}{*}{$1-(144,77,77)$ }\\  
		$[96,48,4]$ &	$2^{29}\cdot 3^2 $&  \\
		$[96,48,8]$ & $(E_{2^3}\times S_3):Z_2$ & $1-(144,77,77)$\\
		$[144,72,8]$ & & \multirow{2}{*}{$1-(144,77,77)$} \\  
		$[96,48,8]$ & $Z_2\times ((S_3\times S_3):Z_2)$ & \\
		$[144,72,12]$ & & \multirow{2}{*}{$1-(144,121,121)$} \\ 
		$[96,48,14]$ & $D_{12}$ &  \\
		$[144,72,12]$ & & $1-(144,121,121)$ \\
		$[144,72,12]$ & & \multirow{2}{*}{$1-(144,89,89$} \\ 
		$[96,48,8]$ & $Z_2\times Z_2\times S_3$ &	  \\
		$[96,48,12]$ & $(S_3\times S_3):Z_2$ & $1-(144,133,133)$ \\
		$[144,72,2]$ & & \multirow{2}{*}{$1-(144,133,133)$} \\ 
		$[96,48,4]$ & $$ &	  \\
		\noalign{\smallskip}\hline
	\end{tabular}
\end{table}

\begin{table}[H]
	\caption{Non-trivial binary pairwise non-equivalent self-orthogonal codes obtained from the orbit matrices of the non-trivial weakly self-orthogonal 1-designs on 165 points using Theorem \ref{tmsl3} (Case 3)}
	\label{sl3t4}
	\begin{tabular}{lll}
		\hline\noalign{\smallskip}
		\textbf{$C$}& \textbf{Aut$C$} or \textbf{$\#$Aut$C$} & \textbf{Design} \\ 
		\noalign{\smallskip}\hline\noalign{\smallskip}
		$[66,33,2]$ &	$2^{64}\cdot 3^{15}\cdot 5^7\cdot 7^4\cdot 11^3\cdot 13^2\cdot 17\cdot 19\cdot 23\cdot 29\cdot 31 $			& \multirow{2}{*}{$1-(165,129,129)$}\\  
		$[30,15,2]$ &  $2^{26}\cdot 3^6\cdot 5^3\cdot 7^2\cdot 11\cdot 13$	 					&  	\\
		$[66,33,4]$ & $E_{2^3}\times((Z_2\times(Z_4\times Z_2):Z_2):Z_2)):Z_2) $ 			& \multirow{2}{*}{$1-(165,81,81)$}\\  
		$[30,15,6]$ & $Z_2 \times (E_{2^4} : A_8) $ &   	\\
		$[66,33,2]$ & $$ & $1-(165,49,49)$ \\
		\noalign{\smallskip}\hline
	\end{tabular}
\end{table}

\begin{table}[H]
	\caption{Non-trivial binary pairwise non-equivalent self-orthogonal codes obtained from the orbit matrices of the non-trivial weakly self-orthogonal 1-designs on 165 points using Theorem \ref{tmsl4} (Case 4)}
	\label{sl4t1}
	\begin{tabular}{lll}
		\hline\noalign{\smallskip}
		\textbf{$C$}& \textbf{Aut$C$} or \textbf{$\#$Aut$C$} & \textbf{Design} \\ 
		\noalign{\smallskip}\hline\noalign{\smallskip}
		$[34,9,6]$ & $E_{2^2}\times((E_{2^2}\times((Z_2\times D_8):Z_2)):Z_2$ & \multirow{2}{*}{$1-(165,109,109)$}\\ 
		$[16,5,8]$* & $E_{2^4}:A_8 $ 	 &  	\\
		$[34,3,10]$ &  $2^{38}\cdot 3^{13}\cdot 5^5\cdot 7^2\cdot 11^2$& $1-(165,61,61)$\\
		$[34,11,4]$ & $E_{2^2}\times ((E_{2^2}\times D_8):Z_2)$ & $1-(165,85,85)$ \\
		\noalign{\smallskip}\hline
	\end{tabular}
\end{table}

\noindent
Using Theorem \ref{fixpoint1} and Theorem \ref{ftsl3}, we constructed binary self-orthogonal codes from the orbit matrices of the non-trivial weakly self-orthogonal 1-designs on less than 165 points (including). In order to construct the orbit matrices we determined orbits of cyclic subgroup $Z_2$ of the group $M_{11}$ acting with orbits of  lengths $1$ and $2.$

\begin{table}[H]
	\caption{Non-trivial binary pairwise non-equivalent self-orthogonal codes obtained from the orbit matrices of the non-trivial self-orthogonal 1-designs on 22 points using Theorem \ref{fixpoint1} (Case 1)}
	\label{fixedpoints22}
	\begin{tabular}{lll}
		\hline\noalign{\smallskip}
		\textbf{$C$}& \textbf{Aut$C$} or \textbf{$\#$Aut$C$} & \textbf{Design} \\ 
		\noalign{\smallskip}\hline\noalign{\smallskip}
		$[6,2,4]$*& $Z_2\times S_4$  & \multirow{2}{*}{$1-(22,10,10)$} \\ 
		$[8,4,2]$  & $Z_2\times S_4$ &  \\
		$[6,3,2]+$  & $Z_2\times S_4$   & $1-(22,2,1)$ \\
		\noalign{\smallskip}\hline
	\end{tabular}
\end{table}

\begin{table}[H]
	\caption{Non-trivial binary pairwise non-equivalent self-orthogonal codes obtained from the orbit matrices of the non-trivial self-orthogonal 1-designs on 66 points using Theorem \ref{fixpoint1} (Case 1)}
	\label{fixedpoints66}
	\begin{tabular}{lll}
		\hline\noalign{\smallskip}
		\textbf{$C$}& \textbf{Aut$C$} or \textbf{$\#$Aut$C$} & \textbf{Design} \\ 
		\noalign{\smallskip}\hline\noalign{\smallskip}
		$[10,2,4]$& $Z_2\times S_4\times S_4$ & \multirow{2}{*}{$1-(66,20,20)$} \\ 
		$[28,4,10]$  & $509607936$ &  \\
		$[10,3,4]+$  & $Z_2\times S_4\times S_4$  & $1-(66,46,46)$ \\
		\noalign{\smallskip}\hline
	\end{tabular}
\end{table}

\begin{table}[H]
	\caption{Non-trivial binary pairwise non-equivalent self-orthogonal codes obtained from the orbit matrices of the non-trivial self-orthogonal 1-designs on 110 points using Theorem \ref{fixpoint1} (Case 1)}
	\label{fixedpoints110}
	\begin{tabular}{lll}
		\hline\noalign{\smallskip}
		\textbf{$C$}& \textbf{Aut$C$} or \textbf{$\#$Aut$C$} & \textbf{Design} \\ 
		\noalign{\smallskip}\hline\noalign{\smallskip}
		$[52,20,4]$ &  $2^{40} \cdot 3^6 $	 					&  $1-(110,72,36)$	\\
		$[52,4,18]$ &  $2^40 \cdot 3^{16} \cdot 5^4 $	 					&  $1-(110,36,18)$	\\
		$[52,24,2]$ &  $2^{49} \cdot 3^{11} \cdot 5^4 \cdot 7^3 \cdot 11^2 \cdot 13 \cdot 17 \cdot 19 \cdot 23 $	& $1-(110,2,1)$ \\
		$[14,2,4]$ &  	$Z_2\times S_8 \times S_4$	 & \multirow{2}{*}{$1-(110,36,18)$}\\ 
		$[48,4,18]$ &   $2^{37} \cdot 3^{15} \cdot 5^4  $	 										&  \\
		$[14,7,2]$ &  	$Z_2\times ((E_{2^6}: A_7):Z_2) $			 & \multirow{2}{*}{$1-(110,2,1)$}\\ 
		$[48,24,2]$ &   $2^{46} \cdot 3^{10} \cdot 5^4 \cdot 7^3 \cdot 11^2 \cdot 13 \cdot 17 \cdot 19 \cdot 23  $	 										&  \\
		$[14,4,4]$ & $(((((Z_2 \times D_8) : Z_2) : Z_3) : Z_2) : Z_2) \times S_6 $		 & \multirow{2}{*}{$1-(110,36,36)$}\\ 
		$[48,20,4]$	&   $ 2^{37} \cdot 3^5$	 & \\
		$[14,5,2]$ & $ (((((Z_2 \times D_8) : Z_2) : Z_3) : Z_2) : Z_2) \times S_6 $			 & $1-(110,39,19)$  \\
		$[14,6,4]+$ & $Z_2\times(((E_{2^6}:A_7):Z_2)$ & $1-(110,72,72)$  \\
		$[14,3,4]$ & $Z_2\times S_8 \times S_4$ &  $1-(110,74,37)$  \\
		$[52,4,10]$ &  $2^{45}\cdot 3^{22} \cdot 5^{10} \cdot 7^5 \cdot 11$		& $1-(110,10,1)$  \\
		$[52,24,2]$ &  $2^{27}\cdot 3^5$	 		& $1-(110,18,9)$\\
		\noalign{\smallskip}\hline
	\end{tabular}
\end{table}

\begin{table}[H]
	\caption{Non-trivial binary pairwise non-equivalent self-orthogonal codes obtained from the orbit matrices of the non-trivial self-orthogonal 1-designs on 132 points using Theorem \ref{fixpoint1} (Case 1)}
	\label{fixedpoints132}
	\begin{tabular}{lll}
		\hline\noalign{\smallskip}
		\textbf{$C$}& \textbf{Aut$C$} or \textbf{$\#$Aut$C$} & \textbf{Design} \\ 
		\noalign{\smallskip}\hline\noalign{\smallskip}
		$[60,24,8]$ &   $Z_2\times S_4$	 					&  $1-(132,20,20)$	\\ 
		$[60,8,6]$ & $2^{49} \cdot 3^{23} \cdot 5^{11}\cdot 7^2 \cdot 11 $	 & $1-(132,6,1)$ 	\\ 
		$[60,28,6]$ &   $S_4$  &  $1-(132,40,40)$	\\ 
		$[12,3,4]$ & $ ((((((A_4 \times A_4) : Z_2) \times A_4) : Z_2) : Z_3) : Z_2) : Z_2 $ 		& \multirow{2}{*}{$1-(132,12,1)$}\\  
		$[60,4,12]$ & $ 2^{53} \cdot 3^{26} \cdot 5^{10} \cdot 7^5 \cdot 11^5 $	 					&  	\\ 
		$[12,3,6]*$ & $((((Z_2 \times (E_{2^4} : Z_2)) : Z_2) : Z_3) : Z_2) : Z_2$ &	 $1-(132,66,6)$	\\ 
		$[60,20,12]$ & $S_4$ 			&  $1-(132,32,32)$	\\ 
		$[12,2,8]*$ & $ ((((((A_4 \times A_4) : Z_2) \times A_4) : Z_2) : Z_3) : Z_2) : Z_2 $  & $1-(132,100,100)$ \\ 
		$[12,4,4]$ &   $(((((Z_2 \times (E_{2^4} : Z_2)) : Z_2) : Z_3) : Z_2) : Z_2) : Z_2$ &  \multirow{2}{*}{$1-(132,32,16)$	}\\  
		$[60,20,8]$ & $ 2^{34} \cdot 3^2 $ &  	\\ 
		$[60,24,4]$ &  $2^{52}\cdot 3^6 $ &  $1-(132,60,30)$	\\ 
		$[60,8,20]$ &  $2^{34} \cdot 3^2 $ &  $1-(132,30,15)$	\\ 
		$[60,28,2]$ & $2^{56}\cdot 3^{14}\cdot 5^6 \cdot 7^4 \cdot 11^2 \cdot 13^2 \cdot 17 \cdot 19 \cdot 23 $ &  $1-(132,2,1)$	\\ 
		$[60,4,20]$ &   $((((((A_4 \times A_4) : Z_2) \times A_4) : Z_2) : Z_3) : Z_2) : Z_2$ &  $1-(132,92,46)$	\\ 
		\noalign{\smallskip}\hline
	\end{tabular}
\end{table}

\begin{table}[H]
	\caption{Non-trivial binary pairwise non-equivalent self-orthogonal codes obtained from the orbit matrices of the non-trivial weakly self-orthogonal 1-designs on 66 points using Theorem \ref{ftsl3} (Case 3)}
	\label{sl3fp66}
	\begin{tabular}{lll}
		\hline\noalign{\smallskip}
		\textbf{$C$}& \textbf{Aut$C$} or \textbf{$\#$Aut$C$} & \textbf{Design} \\ 
		\noalign{\smallskip}\hline\noalign{\smallskip}
		$[20,10,2]$ & $((E_{2^5} : A_6) : Z_2) \times (((((Z_2 \times D_8) : Z_2) : Z_3) : Z_2) : Z_2)$ &  \multirow{2}{*}{$1-(66,21,21)$}\\ 
		$[56,28,4]$ &  $2^{39}\cdot 3^5$	 &  	\\ 
		$[20,10,4]$ & $((((Z_2 \times D_8) : Z_2) : Z_3) : Z_2) \times ((E_{2^5}) : A_6) : Z_2)$ &	$1-(66,45,45)$ \\
		\noalign{\smallskip}\hline
	\end{tabular}
\end{table}

\begin{table}[H]
	\caption{Non-trivial binary self-orthogonal codes of the orbit matrices obtained from the non-trivial weakly self-orthogonal 1-designs on 110 points using Theorem \ref{ftsl3} (Case 3)}
	\label{sl3fp110}
	\begin{tabular}{lll}
		\hline\noalign{\smallskip}
		\textbf{$C$}& \textbf{Aut$C$} or \textbf{$\#$Aut$C$} & \textbf{Design} \\ 
		\noalign{\smallskip}\hline\noalign{\smallskip}
		$[104,52,2]$ & $$ & $1-(110,73,73)$ \\
		$[104,52,2]$ & $ $ 	 & $1-(110,37,37)$ 	\\ 
		$[28,14,2]$ &$2^{24}\cdot 3^4\cdot 5^2 \cdot 7 $& \multirow{2}{*}{$1-(110,37,37)$}\\ 
		$[96,48,4]$ &  $2^{94}\cdot 3^{22}\cdot 5^{10} \cdot 7^6 \cdot 11^4 \cdot 13^3\cdot 17^2\cdot 19^2\cdot 23^2\cdot 29\cdot 31\cdot 37\cdot 41\cdot 43\cdot 47$	 &  	\\
		$[28,14,4]$ &  $2^{23}\cdot 3^4 \cdot 5^2 \cdot 7 $& \multirow{2}{*}{$1-(110,73,73)$}\\ 
		$[96,48,4]$ & $ $ 	 &  	\\
		$[104,52,2]$ & $ $ 	 &  $1-(110,9,9)$	\\
		$[104,52,2]$ & $ $	 &  $1-(110,73,73)$	\\
		$[104,52,4]$ & $ $	 &  $1-(110,81,81)$	\\
		$[104,52,4]$ & $ $ 	 &  $1-(110,29,29)$	\\ 
		$[104,52,2]$ & $ $  &  $1-(110,37,37)$	\\
		$[104,52,2]$ & $ $  &  $1-(110,101,101)$	\\ 
		\noalign{\smallskip}\hline
	\end{tabular}
\end{table}

\begin{table}[H]
	\caption{Non-trivial binary self-orthogonal codes obtained from the orbit matrices of the non-trivial weakly self-orthogonal 1-designs on 132 points using Theorem \ref{ftsl3} (Case 3)}
	\label{sl3fp132}
	\begin{tabular}{lll}
		\hline\noalign{\smallskip}
		\textbf{$C$}& \textbf{Aut$C$} or \textbf{$\#$Aut$C$} & \textbf{Design} \\ 
		\noalign{\smallskip}\hline\noalign{\smallskip}
		$[120,60,2]$ & $$ & $1-(132,5,5)$ \\
		$[120,60,10]$ & $$& $1-(132,21,21)$ \\
		$[120,60,2]$ & $ $ 	 & $1-(132,7,7)$ 	\\ 
		$[120,60,8]$ & $$ & $1-(132,25,25)$ \\
		$[120,60,2]$ & $$ & $1-(132,11,11)$ \\
		$[120,60,10]$ & $$ & $1-(132,27,27)$ \\
		$[64,4,32]+$ & $$ & $1-(132,55,55)$ \\
		$[120,60,8]$ & $$ & $1-(132,31,31)$ \\
		$[120,60,8]$ &	$$& $1-(132,101,101)$ \\ 
		$[120,60,10]$  &	$ $& $1-(132,105,105)$ \\ 
		$[120,60,2]$ & $$ & $1-(132,121,121)$ \\
		$[120,60,8]$ & $$ & $1-(132,107,107)$ \\
		$[120,60,2]$ & $$ & $1-(132,125,125)$ \\
		$[120,60,10]$ & $$ & $1-(132,111,111)$ \\
		$[64,4,12]$ & $$ & $1-(132,121,11)$ \\
		$[120,60,2]$ & $$ & $1-(132,127,127)$ \\
		$[120,60,2]$ & $$ & $1-(132,61,61)$ \\
		$[120,60,2]$ & $$ & $1-(132,31,31)$ \\
		$[120,60,2]$ & $$ & $1-(132,91,91)$ \\
		$[120,60,2]$ & $$ & $1-(132,41,41)$ \\
		$[120,60,2]$ & $$ & $1-(132,101,10)$ \\
		$[120,60,2]$ & $$ & $1-(132,71,71)$ \\
		\noalign{\smallskip}\hline
	\end{tabular}
\end{table}

\end{document}